\documentclass[sn-vancouver,iicol]{sn-jnl}

\usepackage{bm}
\usepackage{subfig}
\usepackage{mathtools}

\jyear{2022}%

\theoremstyle{thmstyleone}%
\theoremstyle{thmstyletwo}%

\theoremstyle{thmstylethree}%

\begin{document}

\title[ ]{Block constrained pressure residual preconditioning for two-phase flow in porous media by mixed hybrid finite elements}


\author[1,2]{\fnm{Stefano} \sur{Nardean}}\email{stefano.nardean@unipd.it}

\author[2]{\fnm{Massimiliano} \sur{Ferronato}}\email{massimiliano.ferronato@unipd.it}

\author*[1]{\fnm{Ahmad} \sur{Abushaikha}}\email{aabushaikha@hbku.edu.qa}

\affil[1]{\orgdiv{Division of Sustainable Development, College of Science and Engineering}, \orgname{Qatar Foundation, Hamad Bin Khalifa University}, \orgaddress{\city{Doha}, \country{Qatar}}}

\affil[2]{\orgdiv{Department of Civil, Environmental and Architectural Engineering}, \orgname{University of Padova}, \orgaddress{\city{Padova}, \country{Italy}}}


%
%
\abstract{This work proposes an original preconditioner that couples the Constrained Pressure Residual (CPR) method with block preconditioning for the efficient solution of the linearized systems of equations arising from fully implicit multiphase flow models. This preconditioner, denoted as Block CPR (BCPR), is specifically designed for Lagrange multipliers-based flow models, such as those generated by Mixed Hybrid Finite Element (MHFE) approximations. An original MHFE-based formulation of the two-phase flow model is taken as a reference for the development of the BCPR preconditioner, in which the set of system unknowns comprises both element and face pressures, in addition to the cell saturations, resulting in a $3\times3$ block-structured Jacobian matrix with a $2\times2$ inner pressure problem. The CPR method is one of the most established techniques for reservoir simulations, but most research focused on solutions for Two-Point Flux Approximation (TPFA)-based discretizations that do not readily extend to our problem formulation. Therefore, we designed a dedicated two-stage strategy, inspired by the CPR algorithm, where a block preconditioner is used for the pressure part with the aim at exploiting the inner $2\times2$ structure. 
The proposed preconditioning framework is tested by an extensive experimentation, comprising both synthetic and realistic applications in Cartesian and non-Cartesian domains.}




\maketitle

\section{Introduction}
The physical processes involving the flow of multiple fluids in porous media are mathematically described by a set of coupled nonlinear Partial Differential Equations (PDEs), which are usually numerically solved. A typical solution approach consists of addressing the inner coupling of the PDEs in a Fully Implicit (FI) manner, while the intrinsic nonlinearity is dealt with by a Newton scheme.
Although this strategy is robust and unconditionally stable, it demands the solution of several large-size linearized systems of equations with the Jacobian matrix at each time step, until satisfactory convergence is attained. Numerical evidence shows that this is the most time- and resource-consuming task in a simulation, typically requiring between 60 and 80\% of the total CPU time~\cite{Magras2001,Hu2013,Esler2021}. While the use of other solvers than iterative Krylov subspace methods is in practice unfeasible, given the size of the systems generated in real-world applications, their performance depends mainly on the robustness and efficiency of the preconditioning strategy being supplied and the interplay that develops with the solver itself. This observation explains the importance of equipping Krylov solvers with efficient off-the-shelf preconditioners or, whenever existing techniques are not robust, developing a customized tool tailored for the application at hand. Either way, achieving good efficiency of the linear solver 
allows reducing the runtime and the consumption of computational resources, while attaining the same solution accuracy.

The Constrained Pressure Residual (CPR) preconditioner~\cite{Wallis1983,Wallis1985} is one of the most established strategies for both academic simulators and industrial software (see, for instance, \citep{Zhou2013,Garipov2018,Cremon2020,Klemetsdal2020,Lie2019,Alvestad2022,Rasmussen2021,Schlumberger2020a,Esler2021,Schlumberger2020,Halliburton2014,Lacroix2001,Singh2018}). 
In its original structure, CPR is a two-stage technique, based on the sequential application of two preconditioners for the pressure subproblem (local stage) and the whole Jacobian matrix (global stage), respectively. This strategy is underpinned by a physics-based intuition, i.e., the different character of the pressure and saturation problems, the former being nearly elliptic and the latter hyperbolic. A consolidated scalable option for the preconditioning of the pressure block is Algebraic MultiGrid (AMG)~\citep{Lacroix2003,Cao2005}, which is very well-suited for elliptic problems, while incomplete factorization with zero fill-in (ILU(0)) is the usual choice for the global stage. Preliminary decoupling of pressure from saturation is also an option to improve the accuracy of the algorithm, with the possible downside of undermining the ellipticity of the original pressure problem~\cite{Gries2014}.

Despite being almost forty-year-old, CPR is still the subject of extensive research (see, for instance, the recent review in~\cite{Nardean2022}). One of the main topics consists of extending the CPR algorithm to general purpose reservoir simulators, e.g., \cite{Roy2020,Cremon2020}.
Roy et al.~\cite{Roy2020}, in particular, developed a two-stage CPR-like algorithm for the nonisothermal dead-oil flow model, denoted as Constrained Pressure-Temperature Residual (CPTR). As the name suggests, temperature is elected as a primary variable together with pressure, so a restricted pressure/temperature problem is addressed in the local stage of the CPTR algorithm. Because of the inner $2 \times 2$ structure of this subproblem, a block preconditioner was developed, instead of relying on a single AMG approximation that can possibly prove ineffective. This is also pointed out in the work by Cremon et al.~\cite{Cremon2020}, where such a preconditioning approach was tested. In this article, which focuses on compositional simulations with thermal effects and reactions, another CPR-like algorithm, denoted as CPTR3, was proposed. It consists of three-stages with two local AMG preconditioners for temperature and pressure alone, followed by a global ILU(0) sweep. The CPTR3 preconditioner has been benchmarked against the standard CPR algorithm, showing promising results.

Another research path fostering the integration of the CPR method and block preconditioners concerns the preconditioning of coupled flow/poromechanics problems addressed in an FI fashion. In this context, CPR can be included as a local preconditioner for the flow part within a global block preconditioning framework (see, for instance, \cite{White2019,T.Camargo2021}).

Although the two-phase isothermal flow problem was the original target of the CPR development back in the 80s', it still represents a challenging bench test 
whenever novel and advanced schemes are used for the discretization of the problem PDEs. In this work, we considered a Mixed Hybrid Finite Element (MHFE)~\citep{Brezzi1991} approximation of Darcy's law, coupled with a standard time-implicit Finite Volume (FV) discretization of the mass balance. Following the approach developed in~\cite{Abushaikha2017}, and further applied in~\cite{Abushaikha2020} and \cite{Li2021a}, in which the continuity of the fluxes across the grid interfaces is strongly imposed to ensure the mass balance, the resulting problem is characterized by a $3 \times 3$ block-structured Jacobian matrix, owing to the simultaneous presence of pressure variables defined on both cells and faces, in addition to the saturation unknowns computed on elements. This peculiar matrix format differs from the usual $2 \times 2$ structure generated by Two-Point-Flux-Approximation (TPFA)-based discretizations.

Applying the classical two-stage CPR algorithm, which was precisely designed for TPFA, to such a problem would require approximating the non-symmetric $2 \times 2$ pressure subproblem with a single off-the-shelf AMG preconditioner, which is often not expected to be a robust option. On the contrary, its internal block structure can be exploited, with the introduction of a block preconditioner as a local stage within a broader CPR-like algorithm. The resulting method is denoted as Block CPR (BCPR). A typical strategy for the preconditioning of the systems of equations stemming from single-phase MHFE-based flow simulation, whose $2\times2$ block structure resembles that of our pressure subproblem, consists of eliminating the element pressure unknowns by performing static condensation (see, for instance,~\cite{Kuznetsov2003,Maryska2000}). Then, the resulting matrix can be preconditioned by means of off-the-shelf preconditioners. However, this strategy cannot be applied to our problem formulation, and more generally to problems with compressible fluids, since, unlike the single-phase flow setting, the element pressure block is not diagonal.

In this work, we propose an original approach for addressing the element and face pressure equations.
Building on our previous work on block preconditioners for the single-phase flow problem~\citep{Nardean2020,Nardean2021}, we developed a tool based on the exploitation of approximate versions of the Jacobian matrix decoupling factors to obtain an inexact version of the Schur complement. Although being fully parallelizable, the approximation of such factors is the most expensive task in our algorithm, however, it does not need to be performed at each nonlinear iteration, but only once at the outset of the simulation. Therefore, the associated cost can be amortized throughout the entire simulation.

This article extends the piece of work first presented at the ECMOR 2022 conference in~\cite{Nardean2022a}.
A breakdown of the paper structure is as follows. The preconditioning methodology is described and tested in four applications after the mathematical and numerical model problem are presented. The paper concludes with a discussion of the findings, as well as some preliminary suggestions for further investigation.

\section{Mathematical problem: two-phase flow in porous media}
\label{sec:mat_num_model}
The two-phase flow in porous media, underpinned by some simplifications of the physics, is the model problem considered in this study. The fluids (oil, $o$, and water, $w$), flowing in a compressible medium under isothermal conditions, are assumed to be immiscible and incompressible. Capillarity is also neglected in the model, resulting in the equality of the fluid pressures, i.e., $p_o=p_w=p$. Pressure in the reservoir, $p$, and on the top of the wells, $p_{\text{bh}}$, along with saturation of the wetting phase (water), $S_w$, are elected as model unknowns, consistently with the natural variables' formulation~\citep{Coats1980}.


\subsection{Governing equations}

The set of PDEs describing the flow of multiple fluids in porous media consists of the mass balance equation
\begin{equation}
    \frac{\partial \phi S_{\alpha}}{\partial t} + \nabla \cdot \bm{v}_{\alpha} = f_{\alpha}, \qquad \alpha=o,w \quad \text{on} \ \Omega \times \mathbb{T}, 
    \label{eq:mass_balance_2ph}\\
\end{equation}
and Darcy's law
\begin{multline}
    \bm{v}_\alpha = - \lambda_\alpha K \nabla (p - \gamma_\alpha z), \\
    \alpha=o,w \quad \text{on} \ \Omega \times \mathbb{T}, 
    \label{eq:darcy_eq_mod_2ph}
\end{multline}
written for each phase $\alpha$ in the space and time domains, $\Omega$ and $\mathbb{T}$, respectively.
In equation~\eqref{eq:mass_balance_2ph}, $\phi$ denotes the porosity of the rock, $S_{\alpha}$ and $\bm{v}_\alpha$ are the saturation and velocity of phase
$\alpha$, $t$ is time, and $f_{\alpha}$ is the source/sink term, while, in equation~\eqref{eq:darcy_eq_mod_2ph}, $K$ defines the permeability tensor, $\gamma_{\alpha}$ the specific weight, and $z$ is depth. Furthermore, $\lambda_{\alpha}=k_{r\alpha}/{\mu_{\alpha}}$ denotes the phase mobility factor, where $k_{r\alpha}$ and $\mu_{\alpha}$ are the relative permeability and phase dynamic viscosity, respectively. The factor $\lambda_{\alpha}$ does not depend on pressure, i.e., $\mu_{\alpha}$ is constant, due to the fluid incompressibility assumption. 
The term $f_{\alpha}$ introduces in the reservoir model the action of wells, which are reproduced using the classical Peaceman formulation~\citep{Peaceman1978,Chen2006}.

Additional relationships need to be provided to mathematically close the problem outlined by equations~\eqref{eq:mass_balance_2ph} and \eqref{eq:darcy_eq_mod_2ph}. Specifically, we relied on the analytical Brooks-Corey's model~\citep{Brooks1964} for the relative permeability and the classical relationship
\begin{align}
    \phi = \phi^0 \left[1 + c_r \left(p - p^0 \right) \right],
    \label{eq:porosity}
\end{align}
which relates porosity to the change in pore pressure experienced by the rock~\citep{Aziz1979}.
Here, $c_r$ is the solid phase compressibility and the superscript $0$ denotes the initial conditions of the relevant quantities. The constraint on the sum of the phase saturations
\begin{align}
    \sum_{\alpha=o,w} S_{\alpha} = 1
    \label{eq:sum_satur}
\end{align}
is also an equation of the system. Ultimately, appropriate initial conditions on pressure and saturation,
\begin{subequations}
\begin{align}
    p\rvert_{t=0} &= p^0 & & \text{in} \ \overline{\Omega},   
    \\
    S_w\rvert_{t=0} &= S_w^0 & & \text{in} \ \overline{\Omega},
\end{align}
\label{eq:ini_cond}
\end{subequations}
along with boundary conditions on pressure and fluid fluxes,
\begin{subequations}
\begin{align}
    p &= \overline{p} & & \text{on} \ \Gamma_p \times \mathbb{T}, 
    \\
    - \lambda_{\alpha} K 
    \nabla \left(p - \gamma_{\alpha}z \right) \cdot \bm{n} &= \overline{{v_n}_{\alpha}}  & & \text{on} \ \Gamma_{\bm{v}} \times \mathbb{T}, 
\end{align}
\label{eq:bound_cond}
\end{subequations}
need to be prescribed to give rise to a well-posed mathematical model. In equations~\eqref{eq:ini_cond} and \eqref{eq:bound_cond}, $\Gamma_p$ and $\Gamma_{\bm{v}}$ are the portions of the boundary, $\Gamma$, where Dirichlet and Neumann conditions are enforced, respectively, with $\Gamma_p \cup \Gamma_{\bm{v}} = \Gamma$, and $\overline{\Omega} = \Omega \cup \Gamma$,
while $\bm{n}$ denotes the outer normal vector to $\Gamma$.

\subsection{Weak formulation}
Let $\mathcal{E}^h$ be the collection of non-overlapping hexahedral cells discretizing the reservoir body with $\mathcal{F}^h$ denoting the set of grid faces. The weak form of Darcy's law~\eqref{eq:darcy_eq_mod_2ph} is constructed by means of the MHFE method on the lowest-order $\mathbb{RT}_0$ space~\citep{Raviart1977}, which requires appending pressure unknowns to the center of gravity of both elements and faces, $p^E$ and $\pi$, respectively. Saturation is, instead, computed only on the centroid of the cells. Specifically, the $\mathbb{P}_0$ space with support on the collection of cells is used to reproduce the element pressure and water saturation fields, whereas another $\mathbb{P}_0$ space, defined on faces, plays a similar role for the interface pressure. The 3-D velocity field is approximated in the $\mathbb{RT}_0$ space, which is spanned by piecewise trilinear vector functions, $\bm{\eta}_i^E(x,y,z)$, defined for every face $i$ and element $E$ with support restricted to the cell itself. Darcy's velocity on $E$ is thus expressed in the form:
\begin{equation}
    \bm{v}^h_{\alpha} \rvert_E = \sum_{\ell \in \partial E} \bm{\eta}_\ell^E(x,y,z) q_{\alpha, \ell}^E, \qquad \alpha=o,w,
    \label{eq:velocity_MHFEM}
\end{equation}
where the weight, $q_{\alpha, \ell}^E$, represents the $\alpha$-phase flux across face $\ell$ and $\partial E$ is the collection of element interfaces.
The face unknowns $\pi$ are the hallmark of the MHFE approach over the original Mixed finite element scheme, whose introduction was intended to guarantee the continuity of the normal components of the interface fluxes, and they operate as Lagrange multipliers.

The MHFE weak form of equation~\eqref{eq:darcy_eq_mod_2ph}, expressed for every element of the grid, describes the fluid fluxes across its interfaces based on the pressure and gravitational gradients experienced within the cell itself as \citep{Abushaikha2017}:
\begin{multline}
    \bm{q}_{\alpha}^E = \lambda^*_{\alpha} \left({B^{E}}\right)^{-1} \left[p^E\bm{1} - \bm{\pi}^E - \gamma_{\alpha} \left(z^E\bm{1} - \bm{\zeta}^E \right) \right], \\
    \alpha=o,w.
    \label{eq:fluxes_MHFEM_2ph}
\end{multline}
Here, the vectors $\bm{q}_{\alpha}^E, \ \bm{\pi}^E, \ \bm{\zeta}^E \in \mathbb{R}^{n_f^E}$ collect the fluxes of phase $\alpha$, face pressures and relevant depths, respectively, with $n_f^E$ being the number of element faces. Furthermore, $\bm{1} \in \mathbb{R}^{n_f^E}$ denotes the vector of ones, $\lambda^*_{\alpha}$ is the diagonal matrix containing 
the fluid mobility of the \emph{upstreaming} element for each face, and $B^{E} \in \mathbb{R}^{n_f^E \times n_f^E}$ is the MHFE elementary matrix. The entries of $B^{E}$ are defined as:
\begin{multline}
    B_{ij}^E = \int_{\Omega^E} \left(\bm{\eta}_i^{E}\right)^{T} \left({K^{E}}\right)^{-1} \bm{\eta}^E_j \ d\Omega, \\
    \forall i,j \in \partial E,
    \label{eq:B_ij_2ph}
\end{multline}
where $\Omega^E$ denotes the cell volume. The expression for the interface fluxes in equation~\eqref{eq:fluxes_MHFEM_2ph} is all in all similar to that generated by a Mimetic finite difference discretization (see, for instance,~\cite{Abushaikha2020,Zhang2021}).

In this section, we limited ourselves to presenting just the main aspects of the MHFE method applied to equation~\eqref{eq:darcy_eq_mod_2ph} in a multiphase setting; the interested reader is referred, for instance, to the works~\cite{Maryska1995,Mose1994,Younes2010,Lie2019,Nardean2021} for the full mathematical derivation and further insights.


On the other hand, a classical FV method in space, coupled with a first-order implicit Euler scheme for the time discretization, is used to build the weak form of the mass balance equation~\eqref{eq:mass_balance_2ph}. This is a rather standard approach resulting in the approximate form:
\begin{multline}
    \Omega^E \frac{\phi \left(p^E_n \right)S_{\alpha,n}^E -\phi \left(p^E_{n-1} \right)S_{\alpha,n-1}^E }{\Delta t_{n-1}} + \sum_{i \in \partial E_r} q_{\alpha, i}^{E,E'} \\
    - \Omega^E f_{\alpha}^E = 0, \qquad \forall E \in \mathcal{E}^h,
\label{eq:discr_mass_bal_1}
\end{multline}
where $n-1$ and $n$ denote the previous and actual time steps, respectively, $\Delta t_{n-1} = t_n-t_{n-1}$ is the time increment, and $\partial E_r$ is the set of faces of $E$ not included in the domain boundary. This implies that the reservoir is isolated from the surrounding rock, which is a frequent assumption in reservoir modeling. Moreover, 
\begin{equation}
    q_{\alpha,i}^{E,E'} = \lambda_{\alpha}^* \frac{\bigl(B^{E'}\bigr)^{-1}_{ii} \Lambda_{\alpha,i}^E - \left(B^{E}\right)^{-1}_{ii} \Lambda_{\alpha,i}^{E'}}{\bigl(B^{E}\bigr)^{-1}_{ii} + \bigl(B^{E'}\bigr)^{-1}_{ii}}
    \label{eq:fluxes_MHFEM_STRONG}
\end{equation}
is the flux across the $i$-th face separating elements $E$ and $E'$,
which results from strongly imposing the continuity of the local fluxes, expressed as in equation~\eqref{eq:fluxes_MHFEM_2ph}, on the two sides of interface $i$, $q_{\alpha,i}^E$ and $q_{\alpha,i}^{E'}$, respectively~\citep{Abushaikha2017}. In equation~\eqref{eq:fluxes_MHFEM_STRONG}, we have
\begin{multline*}
    \Lambda_{\alpha,i}^E = L_{B_i^E} \left(p^E-\gamma_{\alpha} z^E \right) \\
    - \sum_{j \in \partial E \setminus \{i\}} \left(B^{E}\right)^{-1}_{ij} \left(\pi_j^E - \gamma_{\alpha}\zeta_j^E \right) 
\end{multline*}
with $L_{B_i^E} = \sum_{j \in \partial E} \left(B^{E}\right)^{-1}_{ij}$.

\subsubsection{The MHFE-FV system of equations}

The governing equations, introduced in the previous section, define a coupled nonlinear problem, $\bm{R}^n=\bm{0}$, at each time step, which is addressed here by an FI strategy. The system equations can be broken down into three groups: $\bm{R}_{\pi} \in \mathbb{R}^{n_f}$, $\bm{R}_p \in \mathbb{R}^{n_E+n_w}$, and $\bm{R}_s \in \mathbb{R}^{n_E}$, where $n_f$, $n_E$, and $n_w$ are the number of faces, elements, and wells in the model. Note that the time step counter $n$ has been dropped in the vectors above to compact the notation. The first set, $\bm{R}_{\pi}=\bm{0}$, aims at ensuring that the total fluxes remain constant across the grid faces, i.e.,
\begin{equation}
    \sum_{\alpha=o,w} \left( q_{\alpha,i}^E + q_{\alpha,i}^{E'} \right) = 0, \qquad \forall i \in \mathcal{F}^h.
    \label{eq:tot_flux_cont}
\end{equation}
The second set, $\bm{R}_p=\bm{0}$, consists of the so-called Implicit Pressure Explicit Saturation (IMPES)-like \emph{pressure equations}, obtained by summing equations~\eqref{eq:discr_mass_bal_1} over the phase index for each cell, in addition to the well constraints, which are included here due to the expected low number of wells as compared to the number of elements.
Ultimately, the set $\bm{R}_s=\bm{0}$ collects the so-called \emph{saturation equations}, i.e., the $n_E$ discretized mass balance equations~\eqref{eq:discr_mass_bal_1} written for the water phase.
The overall size of the nonlinear system of equations, $\bm{R}^n=[\bm{R}_{\pi}, \ \bm{R}_p, \ \bm{R}_s]^{n,T}=\bm{0}$, is $2n_E+n_f+n_w$, with $n_E$, $n_f$, and $n_w$ pressure unknowns computed on elements, faces, and the wells' top perforation, respectively, and $n_E$ saturations.

The nonlinearity of $\bm{R}^n=\bm{0}$ is addressed by a classical Newton scheme
\begin{equation}
    \bm{x}^{n,(m)} = \bm{x}^{n,(m-1)} + \delta \bm{x},
    \label{eq:sol_upd}
\end{equation}
which is stopped when satisfactory convergence to the exact solution is achieved (see also Section~\ref{sec:num_res} for the specific criterion implemented in the simulator).
Here, $\bm{x} = [\bm{x}_{\pi}, \ \bm{x}_{p}, \ \bm{x}_{s}]^T$ is the full set of unknowns arranged in homogeneous groups for the pressure and saturation, and $m$ is the Newton iteration counter. In equation~\eqref{eq:sol_upd}, the solution update $\delta \bm{x}$ is computed by solving the linearized system of equations:
\begin{multline}
    \mathcal{J}^{n,(m-1)} \delta \bm{x} = - \bm{R}^{n,(m-1)} \qquad
    \Rightarrow \qquad \\
    \begin{bmatrix}
        J_{\pi \pi} & J_{\pi p} & J_{\pi s} \\
        J_{p \pi} & J_{p p} & J_{p s} \\
        J_{s \pi} & J_{s p} & J_{s s}
    \end{bmatrix}^{n,(m-1)}
    \begin{bmatrix}
        \delta \bm{x}_{\pi} \\
        \delta \bm{x}_{p} \\
        \delta \bm{x}_{s}
    \end{bmatrix}
    = -
    \begin{bmatrix}
        \bm{R}_{\pi} \\
        \bm{R}_{p} \\
        \bm{R}_{s}
    \end{bmatrix}^{n,(m-1)},
    \label{eq:linearized_syst_prot}
\end{multline}
where $\mathcal{J}^{n,(m-1)}=\frac{\partial \bm{R}^{n,(m-1)}}{\partial \bm{x}}$ is the Jacobian matrix exhibiting a $3 \times 3$ block structure, and $\bm{R}^{n,(m-1)}$ is the residual vector, both updated at the $(m-1)$-th iteration. As to the size of the diagonal blocks, we have $J_{\pi \pi} \in \mathbb{R}^{n_f \times n_f}$, $J_{p p} \in \mathbb{R}^{(n_E + n_w) \times (n_E + n_w)}$, and $J_{s s} \in \mathbb{R}^{n_E \times n_E}$. Before applying the solution update in equation~\eqref{eq:sol_upd}, $\delta \bm{x}$ is properly postprocessed with the Appleyard chop technique~\citep{Younis2011,Schlumberger2020a}, which provides a sort of guidance and stabilization to the behavior of the nonlinear solver.

\begin{figure}[tb]
    \centering
    \subfloat[Non-zero pattern of $\mathcal{J}$ as per equations ~\eqref{eq:Jacobi_2x2} and \eqref{eq:Jacobi_2x2_1}.]{\qquad \includegraphics[width=0.25\textwidth]{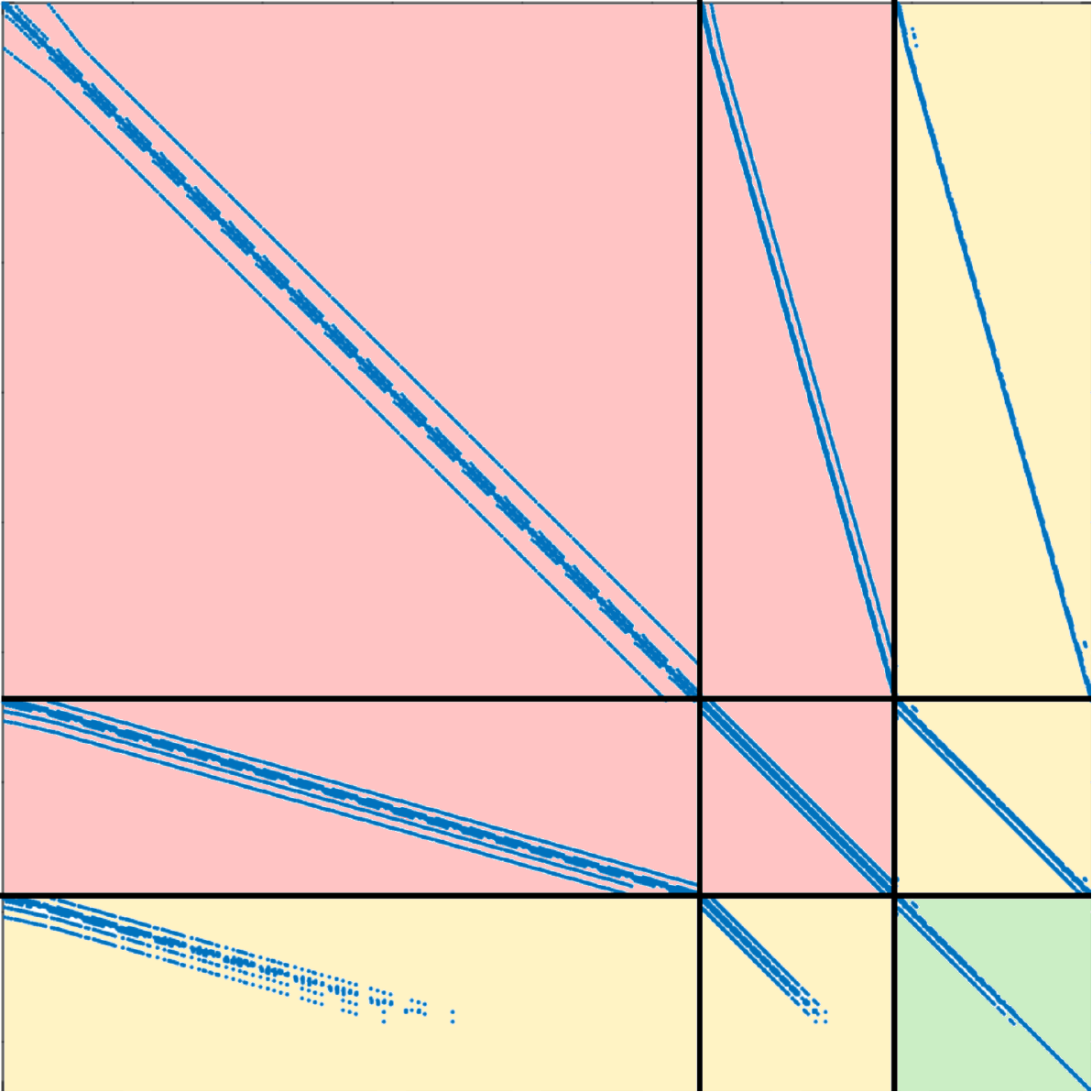} \qquad
    }
    \qquad
    \subfloat[Non-zero pattern of $\mathcal{J}$ after element-wise reordering of $\bm{x}_p$ and $\bm{x}_s$ unknowns.]{\qquad
    \includegraphics[width=0.25\textwidth]{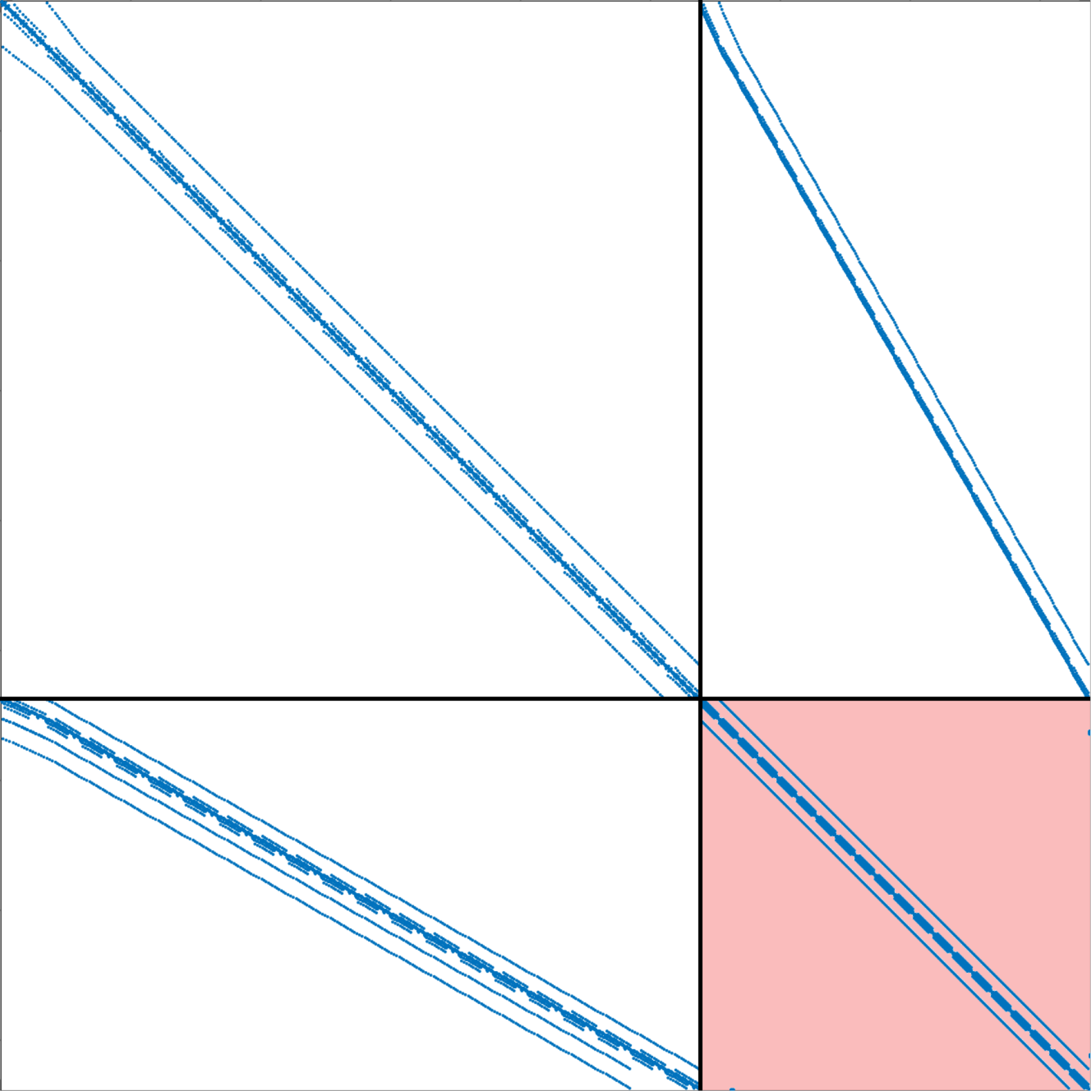} \qquad
    }
    \caption{Non-zero pattern of $\mathcal{J}$ with different unknowns ordering. In panel (b), the block in position (2,2), highlighted in red, roughly corresponds to the system matrix obtained from a standard TPFA-FV discretization with its typical seven-point stencil. The off-diagonal blocks express the coupling between element and face unknowns. Notice that block (2,1) is denser than (1,2) as a result of strongly imposing the flux continuity in the mass balance equations.}
    \label{fig:Jac_example}
\end{figure}

Inspection of the Jacobian matrix reveals that $J_{\pi s}=0$ when gravity is neglected, otherwise all blocks are non-zero. Moreover, under the same condition, $J_{\pi p}$ and $J_{\pi \pi}$ do not depend on $\lambda_{\alpha}$. This means that these blocks remain constant throughout the simulation with $J_{\pi \pi}$ being also Symmetric Positive Definite (SPD), while the whole Jacobian is always non-symmetric, regardless of whether gravitational forces are included in the model or not.
For the design of the BCPR preconditioner, it is convenient to rearrange the Jacobian matrix into a $2 \times 2$ structure by grouping together the pressure blocks, yielding:
\begin{equation}
    \mathcal{J} =
    \begin{bmatrix}
        J_{PP} & J_{Ps} \\
        J_{sP} & J_{ss}
    \end{bmatrix},
    \label{eq:Jacobi_2x2}
\end{equation}
where
\begin{multline}
    J_{PP} =
    \begin{bmatrix}
       J_{\pi \pi} & J_{\pi p}  \\
        J_{p \pi} & J_{p p} 
    \end{bmatrix}
    , \qquad
    J_{P s} = [ J_{\pi s} \ J_{ps}]^T, \\
    J_{sP} = [J_{s \pi} \ J_{s p}].
    \label{eq:Jacobi_2x2_1}
\end{multline}
A sketch of the non-zero pattern of $\mathcal{J}$ is shown in Figure \ref{fig:Jac_example}a.

\section{The Block CPR (BCPR) preconditioner}
The CPR-like family of preconditioners is one of the most successful tools for the efficient preconditioning of linearized systems of equations in reservoir simulations. 
The use of TPFA for the inter-element fluxes approximation, with a single pressure unknown per cell, influenced the structure of the CPR algorithm, where a computationally efficient and scalable choice for the approximation of the elliptic pressure block alone is offered by AMG. 
However, applying the classical CPR algorithm to the linearized system in equation~\eqref{eq:linearized_syst_prot} is not optimal since $J_{PP}$ is a $2 \times 2$ non-symmetric block matrix, compounding both element and face variables, 
for which an approximation by means of existing off-the-shelf AMG functions is often ineffective.
On the contrary, such a block structure can be rather exploited by replacing AMG for the whole pressure part with a block preconditioner relying, in turn, on AMG approximations for the application of the leading block and Schur complement. In essence, we use the classical two-stage CPR multiplicative framework:
\begin{equation}
    \mathcal{M}_{\textup{BCPR}}^{-1} = \mathcal{M}_2^{-1} \left[I - \mathcal{J}\mathcal{M}_1^{-1} \right] + \mathcal{M}_1^{-1},
    \label{eq:BCPR_prec}
\end{equation}
where 
\begin{equation}
    \mathcal{M}_1^{-1} = \text{diag}\left(\mathcal{J}\right)^{-1}
    \label{eq:first_stage_prec}
\end{equation}
and
\begin{equation}
    \mathcal{M}_2^{-1} =
    \begin{bmatrix}
        \mathcal{M}_{PP}^{-1} & \\
                              & 0_s
    \end{bmatrix}
    \label{eq:second_stage_prec}
\end{equation}
are the first- and second-stage preconditioners, respectively. In equation~\eqref{eq:first_stage_prec}, $\text{diag}\left(\mathcal{J} \right)^{-1}$ is the Jacobi preconditioner applied to $\mathcal{J}$, while, in~\eqref{eq:second_stage_prec}, $\mathcal{M}_{PP}^{-1}$ denotes the block preconditioner for the overall pressure subproblem 
and $0_s$ is the zero matrix in the space of saturations.

By inspecting equations~\eqref{eq:BCPR_prec}-\eqref{eq:second_stage_prec},
we notice that two modifications, in addition to the block preconditioner for the pressure subproblem, have been introduced with respect to the most classical CPR framework:
\begin{enumerate}
\item The order of the stages is inverted, i.e., the global stage is carried out before the local one. This is done following the suggestion in~\cite{Roy2020} (and previously applied in \cite{Bui2017}) in order to provide some decoupling of pressure from saturation. The decoupling task, in fact, is not explicitly performed in a preliminary stage, so as to preserve the original algebraic properties and stencil of the $J_{PP}$ block, which is key for building $\mathcal{M}_{PP}^{-1}$; 
\item The global stage is carried out by a simple Jacobi preconditioner instead of the classical ILU(0). 
System~\eqref{eq:linearized_syst_prot} embraces both element and face unknowns, which do not trivially allow for a compact matrix band. 
Numerical results showed that ILU(0), applied to the Jacobian matrix reordered as in Figure~\ref{fig:Jac_example}b, 
does not provide an effective approximation and often may introduce significant round-off errors. 
Although other remedies may include allowing some fill-in in the ILU decomposition or applying algebraic reordering techniques, 
we found that 
the inexpensive Jacobi preconditioner can suffice for a satisfactory performance of the BCPR preconditioner, as we will comment more extensively in Section~\ref{sec:test_1}. 
This also appears to be favorable in view of a massive parallelization of the algorithm.
\end{enumerate}

The starting point for the design of a block preconditioner for $J_{PP}$ is its $\mathcal{LDU}$ decomposition. 
By introducing proper approximations for the application of the inverse of $J_{\pi\pi}$ and the Schur complement $S = J_{pp} - J_{p \pi} J_{\pi \pi}^{-1} J_{\pi p}$, the block
preconditioner for the pressure subproblem takes the form:
\begin{multline}
    \mathcal{M}_{PP}^{-1} = \mathcal{U}^{-1} \mathcal{D}^{-1} \mathcal{L}^{-1}
    = \\
    \begin{bmatrix}
       I_{\pi} & -\widetilde{J}_{\pi \pi}^{-1} J_{\pi p} \\
         & I_p
    \end{bmatrix}
    \begin{bmatrix}
       \widetilde{J}_{\pi \pi}^{-1} & \\
        & \widetilde{S}^{-1}
    \end{bmatrix}
    \begin{bmatrix}
       I_{\pi} &  \\
       - J_{p \pi} \widetilde{J}_{\pi \pi}^{-1} & I_p
    \end{bmatrix},
    \label{eq:block_prec_pres}
\end{multline}
where the superscript $\thicksim$ denotes an approximate term. 
Approximating the second part of the Schur complement is usually troublesome, as the unknown and dense term $J_{\pi \pi}^{-1}$ is involved. Building on the Explicit Decoupling Factor Approximation (EDFA) preconditioner, developed in our previous works \citep{Nardean2020,Nardean2021}, we rewrite the Schur complement as
\begin{equation}
    S = J_{pp} + J_{p \pi}F,
    \label{eq:Schur_recast}
\end{equation}
where
\begin{equation}
    F = - J_{\pi \pi}^{-1} J_{\pi p} 
\end{equation}
is the upper decoupling factor in equation~\eqref{eq:block_prec_pres}. The $F$ term can be computed by solving a series of Multiple Right-Hand Side (MRHS) systems with the columns of $J_{\pi p}$ and $J_{\pi \pi}$ as a matrix:
\begin{equation}
    -J_{\pi \pi} F = J_{\pi p}.
\end{equation}
Obviously, the exact calculation of the dense $F$ factor is computationally intensive, so an approximation is sought. In particular, to preserve a workable sparsity, we compute each column of $F$ only at certain locations belonging to a given non-zero pattern. Focusing on the approximation of the $q$-th column of $F$, with $1 \leq q \leq n_E$, let $R_r^{(q)}$ be the restriction operator to the rows belonging to the non-zero pattern and $\bm{l}^{(q)}$ be the $q$-th vector of the canonical basis. The approximation of the column in the restricted space is found by solving the following (small) system:
\begin{equation}
    - J_{\pi \pi}^{(q)} \widetilde{\bm{f}}^{(q)} = R_r^{(q)} \bm{j}_{\pi p}^{(q)},
    \label{eq:res_system}
\end{equation}
where $J_{\pi \pi}^{(q)} = R_r^{(q)} J_{\pi \pi} \bigl( R_r^{(q)} \bigr)^T$ and $\bm{j}_{\pi p}^{(q)} = J_{\pi p} \bm{l}^{(q)}$.
The inexact factor, $\widetilde{F}$, is obtained by gathering all the $\widetilde{\bm{f}}^{(q)}$ contributions after being prolonged to the original space:
\begin{equation}
    \widetilde{F} = \sum_{q=1}^{n_E} \left( R_r^{(q)} \right)^T \widetilde{\bm{f}}^{(q)} \left( \bm{l}^{(q)} \right)^T.
    \label{eq:approx_F}
\end{equation}
Notice that $J_{\pi p} \in \mathbb{R}^{n_f\times(n_E+n_w)}$ so $F$ must share the same size. However, the well and flux continuity equations are mutually decoupled, hence the last $n_w$ columns in $J_{\pi p}$ are zero and we can avoid solving the relevant homogeneous systems~\eqref{eq:res_system}. This explains the upper limit for $q$ in equation~\eqref{eq:approx_F}.

For the approximation $\widetilde{F}$ to be effective, it is crucial to identify the locations of the most significant entries in each column and collect them to form the non-zero pattern, while preserving an adequate sparsity. To this end, we developed two techniques, denoted as \emph{static} and \emph{dynamic}, respectively, the former being more physics-based and the latter fully algebraic.

\begin{figure*}
    \centering
    \subfloat[Original pattern]{\quad
    \includegraphics[scale=0.9]{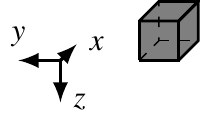} \quad
    } \qquad
    \subfloat[Pattern A]{\quad
    \includegraphics[scale=0.9]{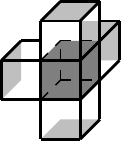} \quad
    } \quad
    \subfloat[Pattern B]{\quad
    \includegraphics[scale=0.9]{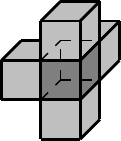} \quad
    } \qquad
    \subfloat[Pattern C]{
    \includegraphics[scale=0.9]{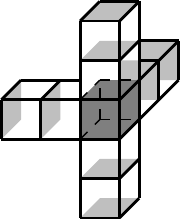}
    } \\
    \subfloat[Pattern D]{
    \includegraphics[scale=0.9]{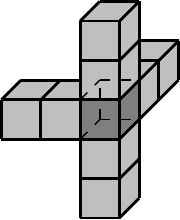}
    } \qquad
    \subfloat[Pattern E]{
    \includegraphics[scale=0.9]{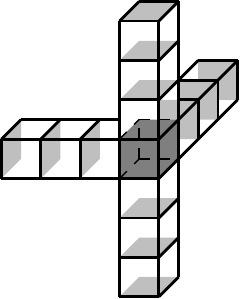}
    } \qquad
    \subfloat[Pattern F]{
    \includegraphics[scale=0.9]{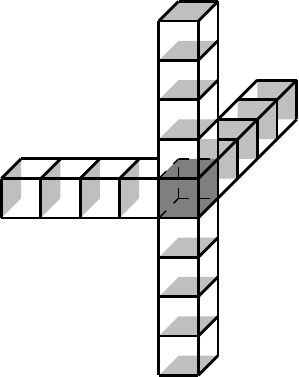}
    }
    \caption{Collection of patterns for the EDFA static approach. Notice that the front and right elements have been removed to improve the picture clarity.}
    \label{fig:static_pat}
\end{figure*}

The static approach is underlain by a graphical interpretation, for which each entry in the $\bm{j}_{\pi p}^{(q)}$ and $\widetilde{\bm{f}}^{(q)}$ vectors is related to a grid face and the non-zeros in $\bm{j}_{\pi p}^{(q)}$ correspond to the faces of element $q$, as shown in Figure~\ref{fig:static_pat}a. Moreover, the problem formalized in equation~\eqref{eq:res_system} is roughly equivalent to solving a flow problem in a portion of the domain whose shape and size is defined by the nonzero pattern. Therefore, starting from $\bm{j}_{\pi p}^{(q)}$'s sparsity pattern, we can customize a broad range of compact patches by including neighboring interfaces such as those displayed in Figure~\ref{fig:static_pat}. The patterns depend on the number of levels of the connection, i.e., how many elements in addition to the central one are involved in each direction (from zero in Figure~\ref{fig:static_pat}a to four in~\ref{fig:static_pat}g), and the presence of the lateral faces, while isotropy in the three directions is preserved. Selecting an appropriate pattern such that the shape and size of the resulting subdomain, where the flow problem is solved, allows for capturing the propagation of the main pressure gradients is key for an accurate approximation of $F$. With the static technique, a prototype pattern is thus supplied beforehand and applied to the whole domain. Numerical tests in~\cite{Nardean2021} showed that this strategy is effective when the grid is Cartesian and the modeler has a robust idea of the possible flow field.

On the other hand, the dynamic technique relies on a recursive algorithm, during which additional column non-zero entries are introduced at the locations where the components of the residual $\bm{r} = \bm{j}_{\pi p}^{(q)} - J_{\pi \pi} (R_r^{(q)})^T \widetilde{\bm{f}}^{(q)}$ are larger in absolute value, starting from the initial non-zero pattern of $\bm{j}_{\pi p}^{(q)}$. Despite being more expensive, the dynamic approach is more flexible and can be used as a black box. The parameters controlling the dynamic process are: (i) $n_{\textup{ent}}$, the number of new entries added to the original pattern, and (ii) $n_{\textup{add}}$, the maximum number of entries introduced at each pattern augmentation step.
%
%
\begin{algorithm}[tb]
\caption{ {\sc Application of the BCPR preconditioner:} $[\bm{v}] = \text{apply\_BCPR}(\tau_i,\mathcal{J},\mathcal{M}_1^{-1},\mathcal{M}_2^{-1},\bm{w})$}
\label{alg:apply_BCPR}
\begin{algorithmic}[1]
\State $\bm{v} = \mathcal{M}_1^{-1} \bm{w}$ \Comment{Application of the first-stage preconditioner}
\State $\bm{r} = \bm{w} - \mathcal{J} \bm{v}$ \Comment{Computing the residual}
\State $[\delta \bm{v}] = \text{apply\_2\_stage}(\tau_i,J_{PP},\widetilde{S},\bm{r})$ \Comment{Application of the second-stage preconditioner to $\bm{r}=[\bm{r}_{\pi},\bm{r}_p,\bm{r}_s]$, $\delta \bm{v} =\mathcal{M}_2^{-1}\bm{r}$}
\State $\bm{v} \gets \bm{v} + \delta \bm{v}$ \Comment{Correction of the first guess}
\end{algorithmic}
\end{algorithm}

\begin{algorithm}[tb]
\caption{ {\sc Application of the second-stage block preconditioner:} $[\delta \bm{v}] = \text{apply\_2\_stage}(\tau_i,J_{PP},\widetilde{S},\bm{r})$}
\label{alg:apply_2_stage}
\begin{algorithmic}[1]
\State $\bm{t}_{\pi} = \text{AMG}(J_{\pi \pi}, \bm{r}_{\pi})$
\State $\bm{t}_p = \bm{r}_p - J_{p \pi} \bm{t}_{\pi}$
\State $\delta \bm{v}_p = \text{solve\_AMG}(\tau_i,\widetilde{S},\bm{t}_p)$
\State $\delta \bm{v} = [\bm{0}, \delta \bm{v}_p, \bm{0}]$
\State $\bm{t}_{\pi} = \bm{r}_{\pi} - J_{\pi p} \delta \bm{v}_p$
\State $\delta \bm{v}_{\pi} = \text{AMG}(J_{\pi \pi}, \bm{t}_{\pi})$
\State $\delta \bm{v} \gets \delta \bm{v} + [\delta \bm{v}_{\pi},\bm{0}, \bm{0}]$
\end{algorithmic}
\end{algorithm}

The action of the BCPR preconditioner to a vector $\bm{w}$ is shown in Algorithm~\ref{alg:apply_BCPR}, while Algorithm~\ref{alg:apply_2_stage} details the application of the block preconditioner at the local stage. AMG, in particular the aggregation-based AGMG presented in~\cite{Napov2012,Notay2010,Notay2012}, is used to approximate block $J_{\pi \pi}$ and the inexact Schur complement
\begin{equation}
    \widetilde{S}=J_{pp}+J_{p\pi}\widetilde{F}.
    \label{eq:Schur_recast_APPROX}
\end{equation}
While a single V-cycle is enough for the application of $\widetilde{J}_{\pi \pi}^{-1}$, an AMG-preconditioned inner Generalized Conjugate Residual (GCR) method~\citep{Eisenstat1983,Jiranek2009} with a loose exit tolerance $\tau_i$ in the range [1.E-5,1.E-4] (achieved in less than 15 iterations) is used for the application of $\widetilde{S}^{-1}$. 
Function solve\_AMG in Algorithm~\ref{alg:apply_2_stage} performs this task.

The choice of using the AGMG function has been made out of convenience since, as opposed to many tools available in the literature, it offers a ready-to-use Matlab interface. In fact, the code implementing the two-phase flow model in compressible media, described in Section~\ref{sec:mat_num_model}, along with the proposed BCPR preconditioner, has been prototyped using Matlab. 
In the future, when transitioning to a low-level programming language implementation of the BCPR preconditioner, more mainstream AMG tools, such as those available in hypre~\cite{Falgout2002} or PETSc~\cite{Balay2022,Balay1997} libraries, will be considered.

In general, the behavior of the global BCPR method depends upon several elements that can be modified or tuned to improve its effectiveness, such as the local- and global-stage preconditioners, the specific AMG tool for the application of the inverse of $J_{\pi \pi}$, as well as the approximation of the Schur complement and of its inverse. In this development stage, we focus on the computation of $\widetilde{S}$, 
setting the other elements as discussed above.


\section{Numerical results}
\label{sec:num_res}
\begin{figure}[tb]
    \centering
    \subfloat[Cartesian grid used in Tests 1-2(a,b) with the horizontal permeability field from the SPE10 data set.]{
    \includegraphics[width=0.47\textwidth]{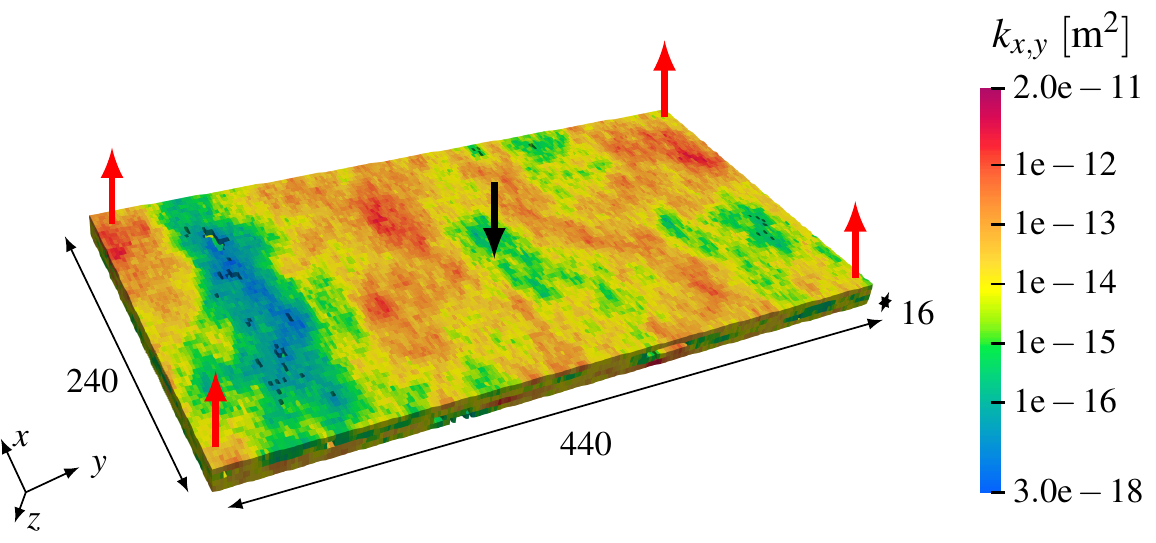}}
    \qquad
    \subfloat[Non-Cartesian grid used in Test 3 with the vertical permability field from the SPE10 data set.]{
    \includegraphics[width=0.47\textwidth]{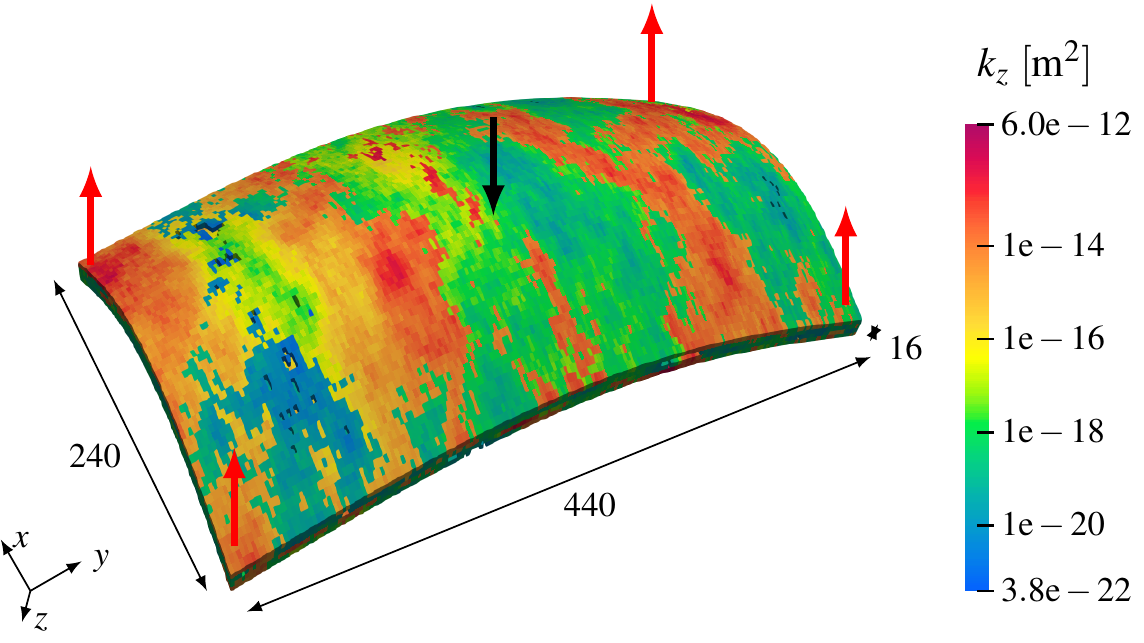}}
    \caption{The two grids used in the tests with a sketch of the simulated scenario. The red and black arrows show the location of production and injection wells, respectively. The horizontal and vertical permeability distributions are superimposed on the grids.}
    \label{fig:Grids}
\end{figure}

The performance of the BCPR preconditioner has been investigated in four test cases, built on the first four layers of the SPE10 data set~\citep{Christie2001} and obtained by varying the original shape of the domain and rock properties.
Figure~\ref{fig:Grids} illustrates the grids, the production scenario and one of the permeability distributions. The domain in Tests 1-2(a,b) is planar and discretized with a Cartesian mesh, while in Test 3 it has been deformed into a dome structure with a non-Cartesian tessellation. The number of elements $n_E$ and faces $n_f$ is the same for both discretizations and equal to 51,741 and 171,070, respectively. The simulated scenario, i.e., the classical five-spot injection/production pattern with an injector in the middle of the reservoir and a producer at each corner, is the same for all tests. The injector operates at a constant rate equal to 20 $\text{m}^3/\text{d}$ and the producers at a constant BHP of 490 and 3,200 kPa for Tests 1, 2a, 3 and 2b, respectively. The wells fully penetrate the reservoir thickness and their radius is equal to 0.1 m. A summary of the tests setup is offered in Table~\ref{tab:test_setup_2_phase}.

The objective of Test 1 is to introduce the BCPR preconditioner and to define a baseline for its setup. To this end, we apply a homogeneous and isotropic permeability distribution ($k_{x,y,z}=\text{1.E-12} \ \text{m}^2$), while the initial porosity is uniform throughout the domain ($\phi^0=\text{2.5E-1}$), and gravity is neglected.
Conversely, Tests 2 and 3 are numerically much more challenging. 
A fully heterogeneous and anisotropic permeability field with nonuniform porosity is introduced in Tests 2a and 2b, which represent a challenging bench test for the performance of the BCPR preconditioner. In test 2b, gravity is also enabled to evaluate the BCPR sensitivity to this feature. 
Finally, in Test 3, a heterogeneous permeability distribution, taken from the SPE10 data set, with an anisotropy ratio up to 3,300,
is applied to the dome-structured domain. The objective of this test is to evaluate the performance of the BCPR preconditioner in a realistic setting and point out the differences in its setup when moving from a Cartesian to a non-Cartesian grid.
\begin{table*}
\begin{minipage}{\textwidth} \small
\caption{Setup of the test cases.}
\label{tab:test_setup_2_phase}
\centering
\begin{tabular}{lccccc}
\toprule
Test & & 1 & 2a & 2b & 3 \\
\midrule
Reservoir type & & Plain & Plain & Plain & Dome \\
\multirow{2}*{Perm. tensor prop.} & & Homogeneous & Heterogeneous & Heterogeneous & Heterogeneous \\
  &  & Isotropic & Anisotropic & Anisotropic & Anisotropic \\
Horiz. perm. range & $\left[ \text{m}^2 \right]$ & 1.E-12 & [3.0E-18, 2.0E-11] & [3.0E-18, 2.0E-11] & [3.0E-18, 2.0E-11] \\
Vert. perm. range & $\left[ \text{m}^2 \right]$ & 1.E-12 & [3.8E-22, 6.0E-12] & [3.8E-22, 6.0E-12] & [1.0E-13, 6.0E-12]\\
Porosity          &                             & 0.25 & [2.6E-05, 5.0E-01] & [2.6E-05, 5.0E-01] & [2.6E-05, 5.0E-01] \\
Oil spec. gravity & $\left[\frac{\text{kPa}}{\text{m}}\right]$ & - & - & 8.00 & - \\
Water spec. gravity & $\left[\frac{\text{kPa}}{\text{m}}\right]$ & - & - & 9.81 & - \\
\botrule
\end{tabular}
\footnotetext{The petrophysical properties of the medium relevant to the two-phase flow model in a compressible porous matrix, common to all the tests, are as follows: rock compressibility $c_r= 5.\text{E-7} \ \text{kPa}^{-1}$, oil dynamic viscosity $\mu_o = 2.3148\text{E-11 kPa d}$, water dynamic viscosity $\mu_w = 1.1574\text{E-11 kPa d}$, irreducible water saturation $S_{wr} = 0$, residual oil saturation $S_{or} = 0$, and initial water saturation $S_{w}^0 = 0$. The relative permeability profiles used in the tests are analytically expressed through Brooks-Corey's model and exhibit a parabolic shape.}
\end{minipage}
\end{table*}

The linearized systems~\eqref{eq:linearized_syst_prot} are solved by means of a right-preconditioned full GMRES~\citep{Saad1986}. Given the relatively low iteration number achieved during the tests, the use of full GMRES appears to be fully warranted. The iterative process is stopped whenever the 2-norm of the relative residuals of the Jacobian system falls below a user-defined threshold, $\tau_l$, i.e., $\| \bm{r}^k \|_2 / \| \bm{r}^0\|_2 < \tau_l$. 
For all tests, $\tau_l = \text{1.E-6}$. The termination criterion for the nonlinear solver relies on the evaluation of both the absolute and relative residual. In particular, given the different nature of the variables, the residual is broken down into three parts, with the iteration process ending when the one of the following conditions is satisfied:
\begin{align*}
    & \max \left\{ \|\bm{R}_{\pi}^{n,(m)}\|_2, \ \|\bm{R}_{p}^{n,(m)}\|_2, \ \|\bm{R}_{s}^{n,(m)}\|_2 \right\} < \tau_{nl,a}, \\
    & \begin{multlined}[b]
        \max \Bigl\{ \|\bm{R}_{\pi}^{n,(m)}\|_2/\|\bm{R}_{\pi}^{n,(0)}\|_2, \ \|\bm{R}_{p}^{n,(m)}\|_2/\|\bm{R}_{p}^{n,(0)}\|_2, \\
        \|\bm{R}_{s}^{n,(m)}\|_2/\|\bm{R}_{s}^{n,(0)}\|_2 \Bigr\} < \tau_{nl,r}.
    \end{multlined}
\end{align*}
For all tests, we set $\tau_{nl,a}=\tau_{nl,r}$ taking the value of $\text{1.E-6}$ for Test 1 and $\text{1.E-5}$ otherwise.
 
In order to provide the most comprehensive picture of the solver computational performance, we identified a set of monitoring parameters defined as follows: 
\begin{enumerate}
\item The number of nonlinear and linear iterations per time step, $N_N$ and $N_l$, respectively; 
\item The total simulation time in seconds per time step, $t_t$, further decomposed into $t_p$ and $t_s$, i.e., the time to build the preconditioner and for GMRES to iterate to convergence, respectively; 
\item The ratio $R_S$ of the number of non-zeros in the Schur complement (equation~\eqref{eq:Schur_recast_APPROX}) computed with the selected pattern for $\widetilde{F}$ with respect to the sparsest one, which corresponds to the nonzero patch of $J_{\pi p}$ (also denoted as original in the sequel), i.e., $\mathtt{nnz}(\widetilde{S})/\mathtt{nnz}(\widetilde{S}_{\textup{orig}})$. 
\end{enumerate}
While the indicators defined in points $1.$ and $2.$ concern the computational efficiency, $R_S$ in point $3.$ provides an indication as to the memory footprint of the preconditioner.
With the aim at measuring both the local and global solver performance, we also introduce the parameters denoted by $\hat{(\cdot)}$ and $\overline{(\cdot)}$, which refer to the cumulative performance throughout the whole simulation and the average performance in a single system, respectively. In particular, the latter refers to the first linearized systems per time step (denoted by the additional subscript 1).
Due to the prototypical Matlab implementation of the code, the results should be regarded as preliminary for the CPU time. An improved absolute CPU time performance is expected with low-level programming languages, such as C++.

\begin{figure}
    \centering
    \subfloat[Systems with $J_{\pi \pi}$ as matrix.]{
    \includegraphics[scale=0.9]{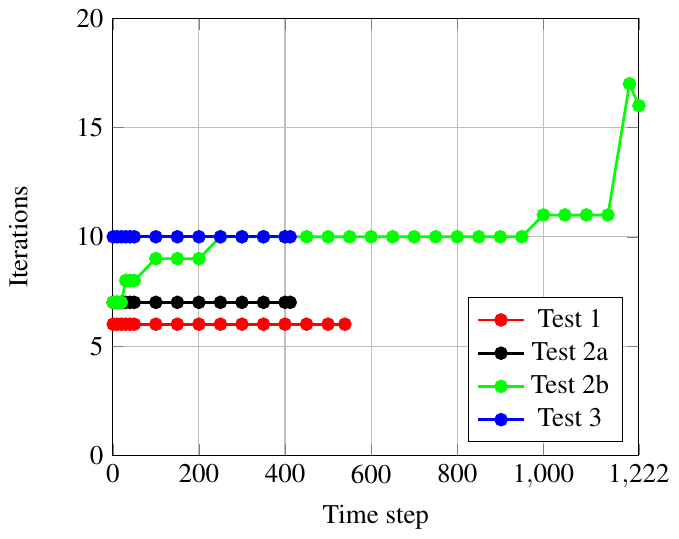}
    } \\
    \subfloat[Systems with $\widetilde{S}$ as matrix.]{
    \includegraphics[scale=0.9]{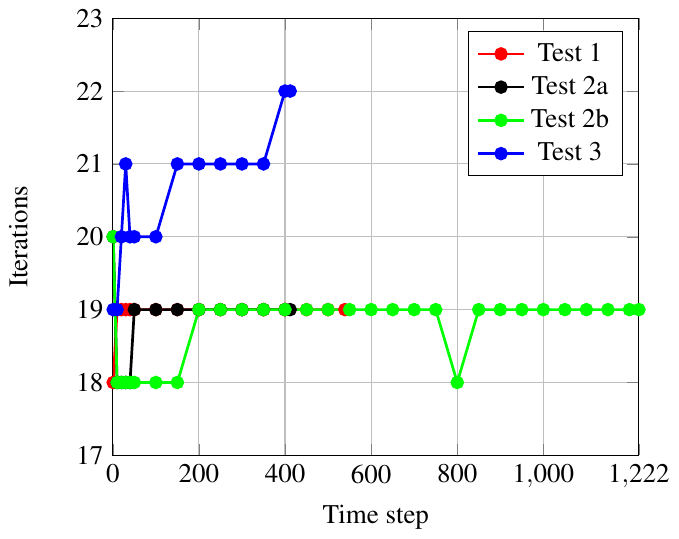}
    }
    \caption{Number of iterations to converge for the systems of equations associated with blocks $J_{\pi \pi}$, (a), and $\widetilde{S}$, (b), generated during the simulations in Tests 1-3. The systems are solved, to a tolerance equal to 1.E-8, by means of GCR preconditioned with AGMG. This solver is also used for SPD systems, in place of the CG method, for consistency. The Schur complement $\widetilde{S}$ is built using pattern A (Figure~\ref{fig:static_pat}b).}
    \label{fig:iter_AGMG_blocks}
\end{figure}

The applicability and effectiveness of the AMG tool to reproduce the action of the inverse of $J_{\pi \pi}$ and $\widetilde{S}$ is a crucial factor. To this aim, we investigate beforehand the performance in the solution of linear systems associated with $J_{\pi \pi}$ and $\widetilde{S}$, extracted during the simulation, for some RHS, using the chosen AMG function as preconditioner. The number of iterations to converge at a given tolerance (equal to 1.E-8 in this analysis) is used to measure the quality of the AMG preconditioner for these local problems. This analysis is carried out at different time steps to evaluate the possible dependency on the problem evolution.
Figure~\ref{fig:iter_AGMG_blocks} summarizes the outcome of this preliminary analysis, where AGMG is used as a preconditioner for the GCR solver~\cite{Eisenstat1983,Jiranek2009}. For the sake of consistency, GCR is used in place of CG even for SPD matrices. Block $J_{\pi \pi} \in \mathbb{R}^{n_f \times n_f}$ arises from an elliptic contribution, so it is well-suited to AMG, and the solver converges in fewer iterations than $\widetilde{S} \in \mathbb{R}^{(n_E + n_w) \times (n_E + n_w)}$. When gravity is not considered in the model (Tests 1, 2a, and 3), $J_{\pi \pi}$ does not depend on mobility, it is constant during the simulation and is also SPD. Otherwise, it is non-symmetric and evolves over time with gravity (Test 2b), and this makes such a scenario computationally more challenging for AMG. As the simulation proceeds, the number of iterations increases slightly and slowly, nonetheless AGMG behaves satisfactorily well. 
As to Figure~\ref{fig:iter_AGMG_blocks}b, we can say that AGMG serves quite satisfactorily as a preconditioner for the Schur complement as well. In particular, the combination of a non-Cartesian grid and heterogeneous anisotropic permeability produces the most challenging scenario. Ultimately, we also observe that the AMG approximation is almost insensitive to the time step size both for $J_{\pi \pi}$ and $\widetilde{S}$. For all the tests, in fact, $\Delta t$ grows, reaching $\Delta t_{\max}$ within the first 40/50 time steps, and, in this range, the iteration count remained almost constant.

While the AMG setup for $\widetilde{S}$ needs to be performed at each nonlinear step, this task can be carried out only once at the outset of the simulation for $J_{\pi \pi}$ when gravity is neglected, since this block is constant during the simulation. On the contrary, when gravity is considered in the model, the AMG setup is required at each nonlinear iteration for both $J_{\pi \pi}$ and $\widetilde{S}$.

The tests are carried out on a workstation equipped with an AMD Ryzen 9 3950X 16-Core processor at 3.49 GHz and 64 GB of RAM.

\subsection{Test 1: Planar reservoir with homogeneous and isotropic permeability}
\label{sec:test_1}
\begin{figure}
    \centering
    \subfloat[Water saturation]{
    \includegraphics[width=0.47\textwidth]{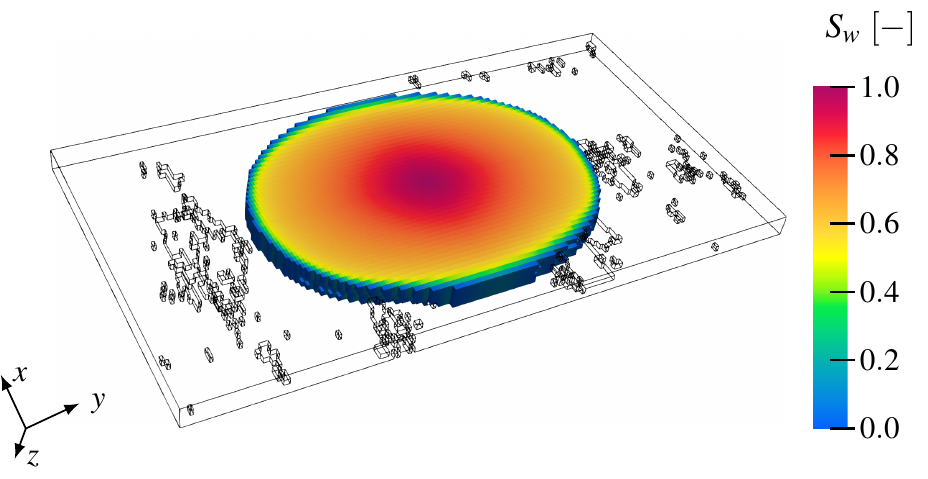}
    }
    \\
    \subfloat[Water velocity]{
    \includegraphics[width=0.47\textwidth]{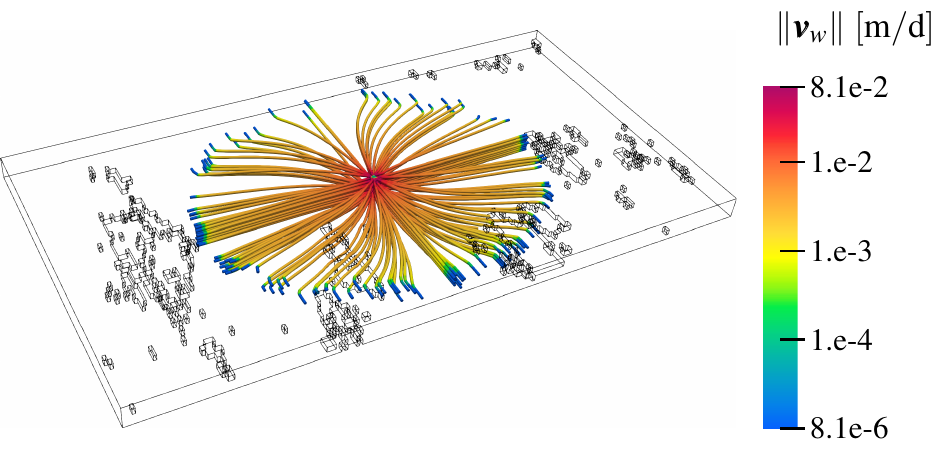}
    }
    \caption{Test 1: Some model insights at the end of the simulation ($t=5,000$ d).}
    \label{fig:T1_insights_2_phase}
\end{figure}
Test 1 reproduces a textbook waterflooding application with production/injection lasting for 5,000 days ($\approx 13.7$ years), which are covered in 539 time steps with a maximum time increment, $\Delta t_{\max}$, equal to 10 days. The Courant–Friedrichs–Lewy (CFL) number is up to 6.25. Figure~\ref{fig:T1_insights_2_phase} offers some model insights in terms of water saturation and water velocity at the end of the simulation. In this application, but also in Tests 2a and 2b, we study the performance of the BCPR preconditioner when the Schur complement is built using the static approach only. In our previous study~\citep{Nardean2021}, in fact, the more expensive dynamic variant proved not necessary when the flow field is moderately regular and the mesh is Cartesian. 
Table~\ref{tab:T1_static} summarizes the main outcome of seven runs performed with the patterns shown in Figure~\ref{fig:static_pat}. 
Not surprisingly, as the pattern is enlarged (runs from 1 to 7), the cumulative time to build the preconditioner, $\hat{t}_p$, increases, whereas the number of linear iterations, $\hat{N}_l$, decreases. The most significant reduction, around 20\%, is from run 1 to 2, i.e., moving from $\widetilde{F}$ built with the original pattern to pattern A. Further patch enlargements give negligible improvements. Actually, the density of the Schur complement (see the $R_S$ column) resulting from the pattern expansion increases the computation and application cost of the preconditioner, so that the best performance is achieved with patterns A-D. This result shows that there is no need for a significant filling of $\widetilde{F}$, with one or two additional levels of the connection, with respect to the original pattern, already enough. Focusing on each single system, convergence is achieved on average in less than 15 iterations. Moreover, we observe that patterns B and D deliver the same results as their counterparts A and C, meaning that they produce the same $\widetilde{F}$ factor. This is due to the block diagonal structure of the $B^E$ matrices (equation~\eqref{eq:B_ij_2ph}) for Cartesian discretizations. Therefore, in that context, the approximation in computing $\widetilde{F}$ is only controlled by the level of the connection and not by the presence of the lateral faces.
\begin{table*}[tb]
    \caption{Test 1: Numerical performance of the static technique.}
    \centering
    \begin{tabular}{cccccccccccc}
    \toprule
    \# & Pat & $R_S$ & $\hat{N}_N$ & $\hat{N}_l$ & $\hat{t}_p$ [s] & $\hat{t}_s$ [s] & $\hat{t}_t$ [s] & $\overline{N}_{l1}$ & $\overline{t}_{p1}$ [s] & $\overline{t}_{s1}$ [s] & $\overline{t}_{t1}$ [s] \\
    \midrule
    1  & Orig & 1   & 1,221 & 22,452 & 68.7  & 2,556.7 & 2,625.4 & 17.7 & 0.1 & 2.0 & 2.1 \\
    2  &   A  & 1.4 & 1,221 & 17,892 & 92.2  & 1,917.7 & 2,009.9 & 14.2 & 0.1 & 1.5 & 1.6 \\
    3  &   B  & 1.4 & 1,221 & 17,892 & 91.4  & 1,911.6 & 2,003.0 & 14.2 & 0.1 & 1.5 & 1.6 \\
    4  &   C  & 1.7 & 1,221 & 17,396 & 104.4 & 1,909.2 & 2,013.6 & 13.7 & 0.1 & 1.5 & 1.6 \\
    5  &   D  & 1.7 & 1,221 & 17,396 & 105.8 & 1,910.8 & 2,016.6 & 13.7 & 0.1 & 1.5 & 1.6 \\
    6  &   E  & 2.0 & 1,221 & 17,341 & 141.7 & 2,163.1 & 2,304.8 & 13.6 & 0.1 & 1.7 & 1.8 \\
    7  &   F  & 2.4 & 1,221 & 17,326 & 154.5 & 2,244.7 & 2,399.2 & 13.6 & 0.1 & 1.7 & 1.8 \\
    \bottomrule
    \end{tabular}
    \label{tab:T1_static}
\end{table*}

\begin{figure}
    \centering
    \subfloat[Linear iterations vs.\ time step]{
    \includegraphics[scale=0.9]{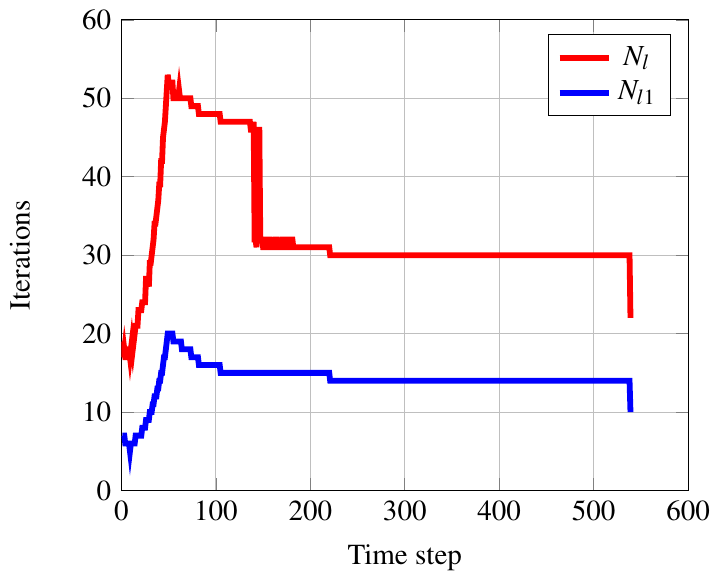}
    } \\
    \subfloat[CPU time vs.\ time step]{
    \includegraphics[scale=0.9]{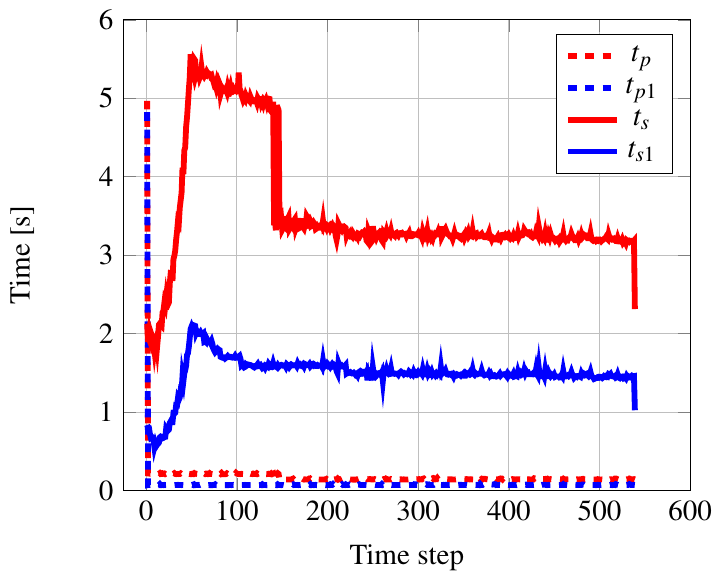}
    }
    \caption{Test 1: Number of linear iterations (a) and CPU time (b) vs.\ the time step index during the full run 2 of Table~\ref{tab:T1_static}.}
    \label{fig:T1_insights}
\end{figure}

The overall performance of the linear solver during the whole simulation is shown in Figure~\ref{fig:T1_insights} in terms of linear iterations and CPU time per time step. The number of nonlinear iterations is 3 until the 140th time step and 2 later on, thus explaining the sudden jump observed in the red profiles. The linear solver displays good performance throughout the full simulation, stabilizing after $\Delta t_{\max}$ is reached, with no impact on the nonlinear convergence, as confirmed in column $\hat{N}_N$ in Table~\ref{tab:T1_static}. Figure~\ref{fig:T1_insights}b also shows, from the peak in the $t_p$ and $t_{p1}$ profiles, that the impact of the computation of $\widetilde{F}$ is at the outset of the simulation. However, such a cost can be effectively amortized during a full run. 
%

\begin{figure}
    \centering
    \includegraphics[scale=0.9]{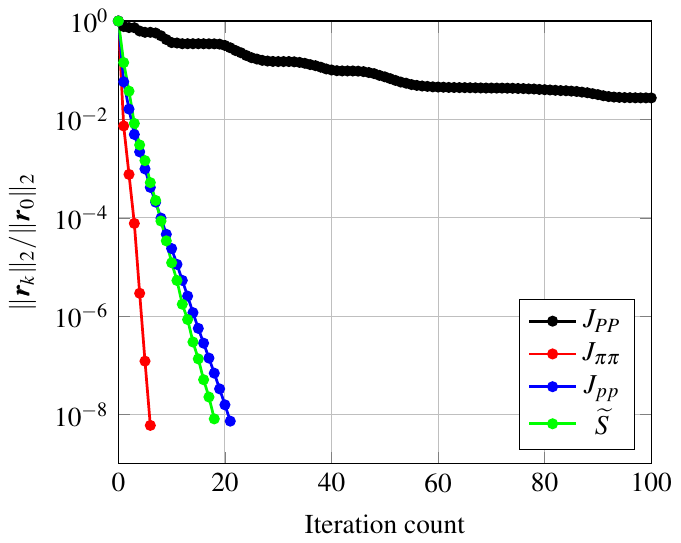}
    \caption{Test 1: Convergence profiles of GCR preconditioned with AGMG for the global pressure problem, $J_{PP}$, leading terms, $J_{\pi \pi}$ and $J_{pp}$, and approximate Schur complement $\widetilde{S}$. The exit tolerance is 1.E-8.}
    \label{fig:T1_conv_profiles}
\end{figure}

The motivation behind the introduction of the BCPR algorithm is the lack of efficiency of conventional general-purpose AMG tools to approximate the whole $2 \times 2$ pressure part of system~\eqref{eq:linearized_syst_prot}. In order to investigate this problem, the system with $J_{PP}$ is solved using the GCR method preconditioned with AGMG for some RHS. The convergence profile is shown in Figure~\ref{fig:T1_conv_profiles} along with those of the diagonal blocks, $J_{\pi \pi}$ and $J_{pp}$, and Schur complement $\widetilde{S}$ for comparison. GCR is also used when the system matrix is SPD ($J_{\pi \pi}$) for consistency. While convergence is fast for $J_{\pi \pi}$, $J_{pp}$, and $\widetilde{S}$, GCR fails to converge within 100 iterations for the system with $J_{PP}$. Therefore, extending the classical CPR method to our model problem~\eqref{eq:linearized_syst_prot} is expected to deliver poor results.

\begin{figure}
    \centering
    \includegraphics[scale=0.9]{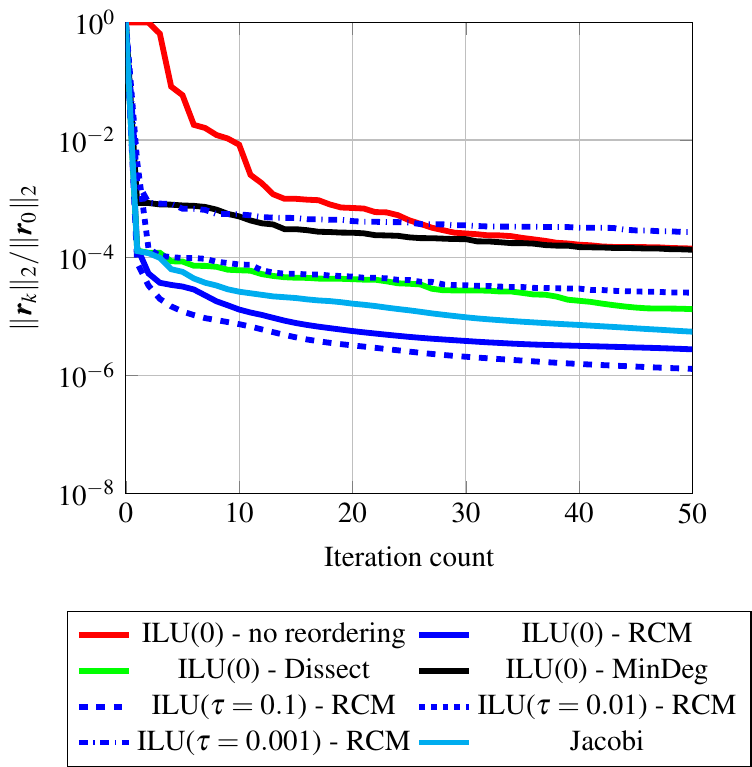}
    \caption{Test 1: Convergence profiles of GMRES preconditioned with Jacobi preconditioner and some variants of ILU factorizations for the Jacobian matrix $\mathcal{J}$. The exit tolerance is 1.E-8.}
    \label{fig:T1_ILU_anal}
\end{figure}

Substituting the ILU(0) factorization with the less expensive Jacobi preconditioner at the global stage is the second main modification, in our BCPR preconditioner, to the original CPR algorithm, since we observed a substantial degradation in the solver convergence as the size of the model is enlarged. The numerical results in Figure~\ref{fig:T1_ILU_anal} help us to support this choice. In this analysis, we solve the system with the Jacobian matrix $\mathcal{J}$ using GMRES preconditioned with some variants of incomplete factorization and the Jacobi preconditioner. The exit tolerance and the maximum number of iterations are set equal to 1.E-8 and 100, respectively. We do not expect GMRES to converge due to the complex block structure of the problem, rather we are interested in the solver behavior during the very first iterations since the global-stage preconditioner is applied only once in the BCPR algorithm. Computing a factorization with zero fill-in directly on the Jacobian matrix ordered as in Figure~\ref{fig:Jac_example}b produces the worst results in terms of the initial convergence rate. 
Applying some algebraic reordering techniques, such as Minimum degree ordering (MinDeg), Nested dissection (Dissect), and Reverse Cuthill-McKee (RCM)~\cite{Saad2003} may help to improve the initial solver behavior. A threshold-based variant of ILU, computed on the RCM-reordered Jacobian matrix, does not allow for a significantly better performance. On the contrary, according to the threshold $\tau$, convergence can degrade, although the factors become very dense. 
On the other hand, the convergence profile with the Jacobi preconditioner (the cyan line) is overall comparable to the best outcome of the reordered ILU factorizations, while being computationally much cheaper 
and embarassingly parallel. 

\subsection{Tests 2a: Planar reservoir with heterogeneous and anisotropic permeability}
\label{sec:test_2a}
\begin{figure*}[tb]
    \centering
    \subfloat[Nonlinear iterations vs.\ time step]{
    \includegraphics[scale=0.9]{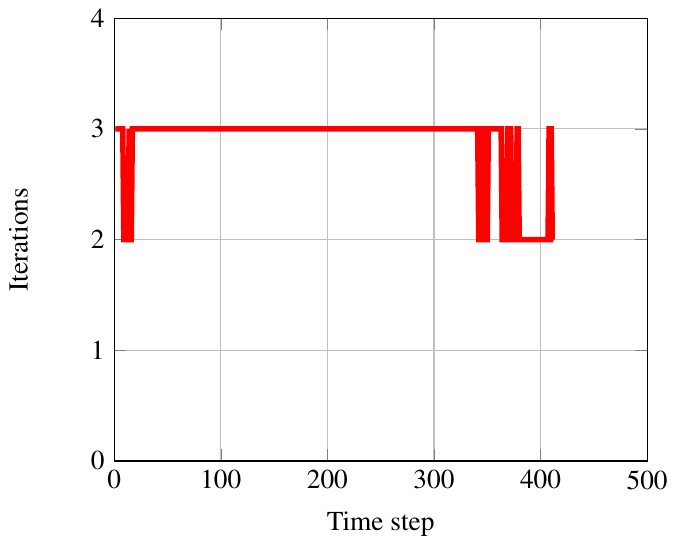}
    } \quad
    \subfloat[Linear iterations vs.\ time step]{
    \includegraphics[scale=0.9]{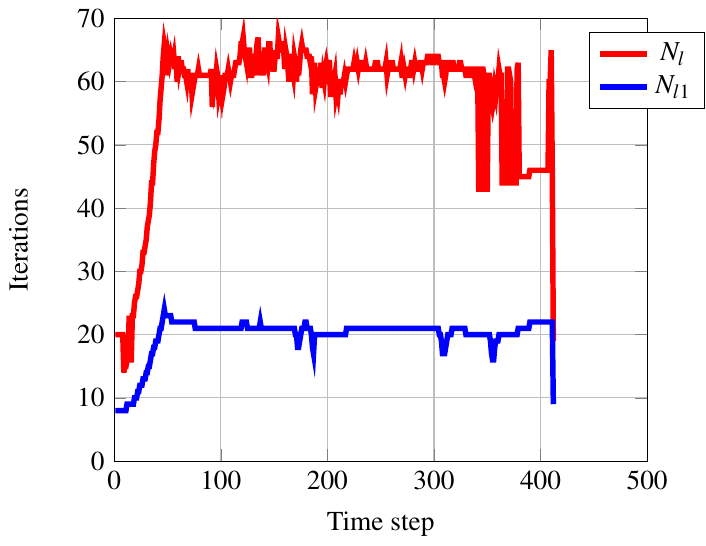}
    } \\
    \subfloat[CPU time vs.\ time step]{
    \includegraphics[scale=0.9]{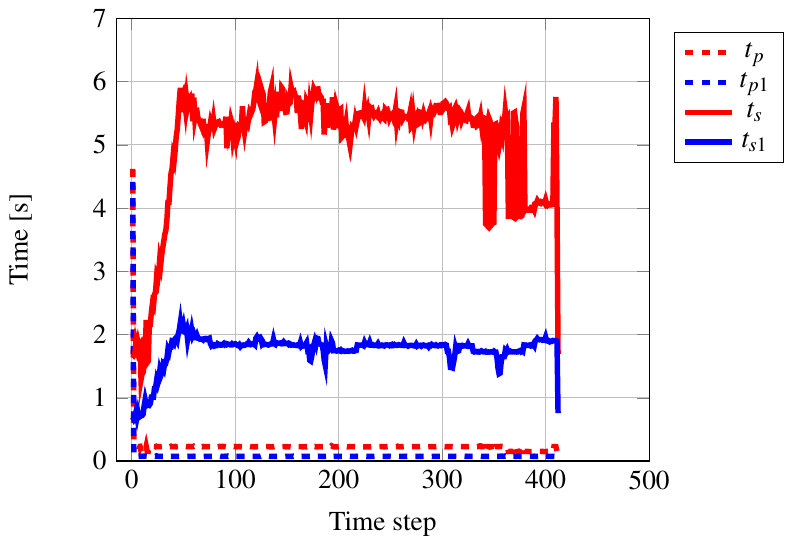}
    }
    \caption{Test 2a: Number of nonlinear (a) and linear (b) iterations and CPU time (c) vs.\ the time step index during the full run 2 of Table~\ref{tab:T2a_static}.}
    \label{fig:T2a_sim_insights}
\end{figure*}

In this test, we consider a challenging scenario characterized by the permeability and porosity distributions of the SPE10 data set. The simulated time is 1,500 days ($\approx 4.1$ years) of continuous production/injection spanned in 412 time steps. The maximum time increment is set equal to 4 days, and the CFL number reaches a value of 19.9.
\begin{figure}
    \centering
    \subfloat[Water saturation]{
    \includegraphics[width=0.47\textwidth]{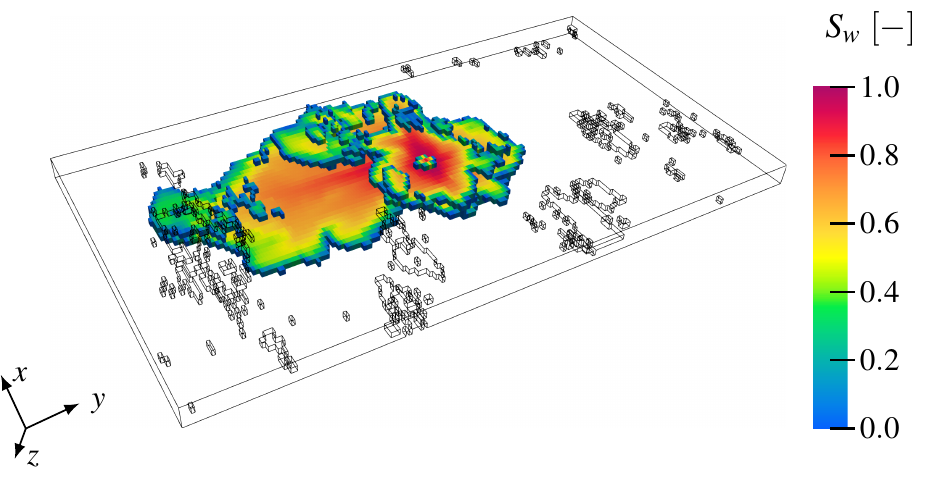}
    }
    \\
    \subfloat[Water velocity]{
    \includegraphics[width=0.47\textwidth]{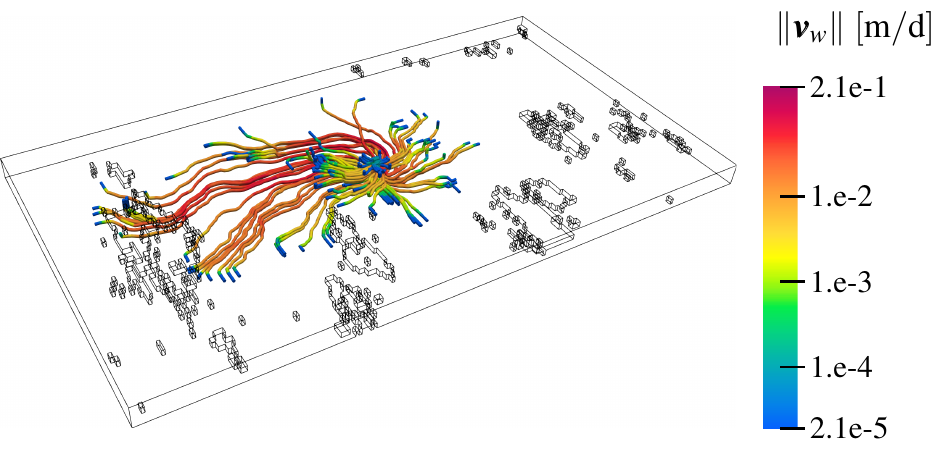}
    }
    \caption{Test 2a: Some model insights at the end of the simulation ($t=1,500$ d).}
    \label{fig:T2a_insights_2_phase}
\end{figure}
Figure~\ref{fig:T2a_insights_2_phase} illustrates the physical outcome of the model at the end of the simulation. The analysis of the BCPR performance with the five patterns without lateral faces is provided in Table~\ref{tab:T2a_static}, where we see that a failure is reported when the Schur complement is built with the original patch. This pattern, in fact, produces a rough approximation of $\widetilde{S}$, which soon loses the algebraic properties that make it suitable for AGMG. Such an issue is easily solved by expanding the element-to-face connection and, once again, pattern A (run 2) delivers the best results in terms of CPU time. The performance of the nonlinear and linear solvers with pattern A for the EDFA method is displayed in Figure~\ref{fig:T2a_sim_insights}, where we can observe a similar behavior as in Test 1.
\begin{table*}[tb]
    \caption{Test 2a: Numerical performance of the static technique.}
    \centering
    \begin{tabular}{cccccccccccc}
    \toprule
    \# & Pat & $R_S$ & $\hat{N}_N$ & $\hat{N}_l$ & $\hat{t}_p$ [s] & $\hat{t}_s$ [s] & $\hat{t}_t$ [s] & $\overline{N}_{l1}$ & $\overline{t}_{p1}$ [s] & $\overline{t}_{s1}$ [s] & $\overline{t}_{t1}$ [s] \\
    \midrule
    1  & Orig & 1   & NC    & -      & -     & -       & -       & -    & -   & -   & -   \\
    2  &  A   & 1.4 & 1,184 & 23,402 & 94.2  & 2,070.9 & 2,165.1 & 19.9 & 0.1 & 1.7 & 1.8 \\
    3  &  C   & 1.7 & 1,183 & 22,488 & 112.6 & 2,209.2 & 2,321.8 & 19.3 & 0.1 & 1.9 & 2.0 \\
    4  &  E   & 2.0 & 1,221 & 23,037 & 138.2 & 2,352.2 & 2,490.4 & 19.3 & 0.1 & 2.0 & 2.1 \\
    5  &  F   & 2.4 & 1,212 & 22,810 & 167.1 & 2,514.8 & 2,681.9 & 19.3 & 0.1 & 2.1 & 2.2 \\
    \bottomrule
    \end{tabular}
    \label{tab:T2a_static}
\end{table*}

\begin{figure}
    \centering
    \subfloat[Cumulative linear iterations during the simulation]{
    \includegraphics[scale=0.9]{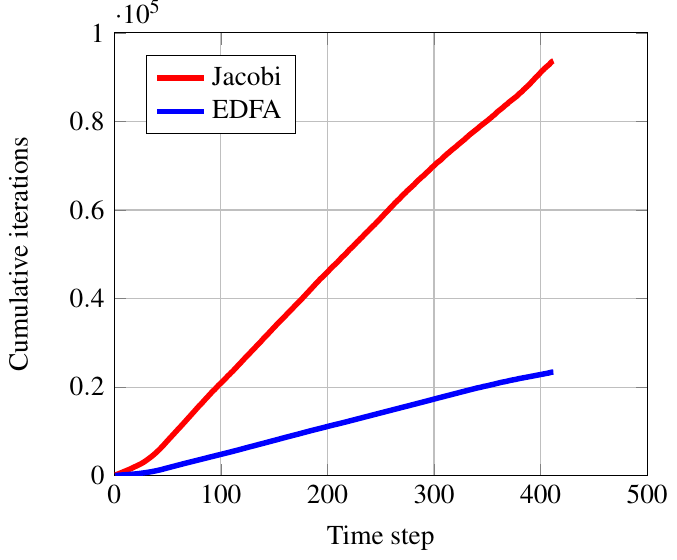}
    } \\
    \subfloat[Cumulative total CPU time during the simulation]{
    \includegraphics[scale=0.9]{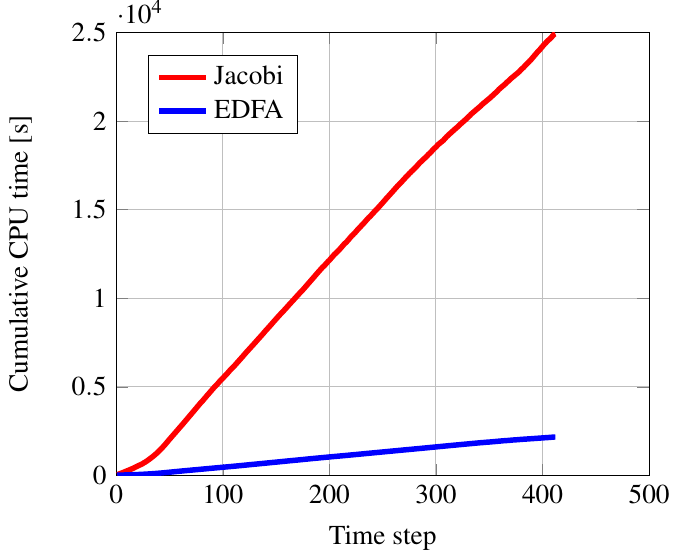}
    }
    \caption{Test 2a: Comparing the performance of different approximations of the Schur complement, using Jacobi preconditioner for $J_{\pi \pi}$ and the best setting for the EDFA method (run 2 in Table~\ref{tab:T2a_static}).}
    \label{fig:T2a_Jacobi_vs_EDFA}
\end{figure}

The computational challenge offered by the SPE10 data set allows us to benchmark, in a realistic setting, our preconditioning solution against a standard approach to approximate the Schur complement.
Performing this task in an effective but inexpensive way is typically one of the most challenging tasks in the design of a block preconditioner. Before resorting to more advanced approaches, Jacobi preconditioner for the inverse of the block in position (1,1) is definitely one of the first attempts, yielding the approximation $\widetilde{S} = J_{p p} - J_{p \pi} \text{diag}(J_{\pi \pi})^{-1} J_{\pi p}$. We tested this option in the BCPR framework and compared the solver performance, in terms of linear iterations and total CPU time, with that obtained with the best EDFA setting (run 2 in Table~\ref{tab:T2a_static}). The EDFA method clearly outperforms the Jacobi-based approach, as shown in Figure~\ref{fig:T2a_Jacobi_vs_EDFA}, which produces a too coarse approximation of the second part of the Schur complement.

\subsection{Tests 2b: Planar reservoir with heterogeneous and anisotropic permeability and gravity}
\label{sec:test_2b}
\begin{figure}
    \centering
    \subfloat[Water saturation]{
    \includegraphics[width=0.47\textwidth]{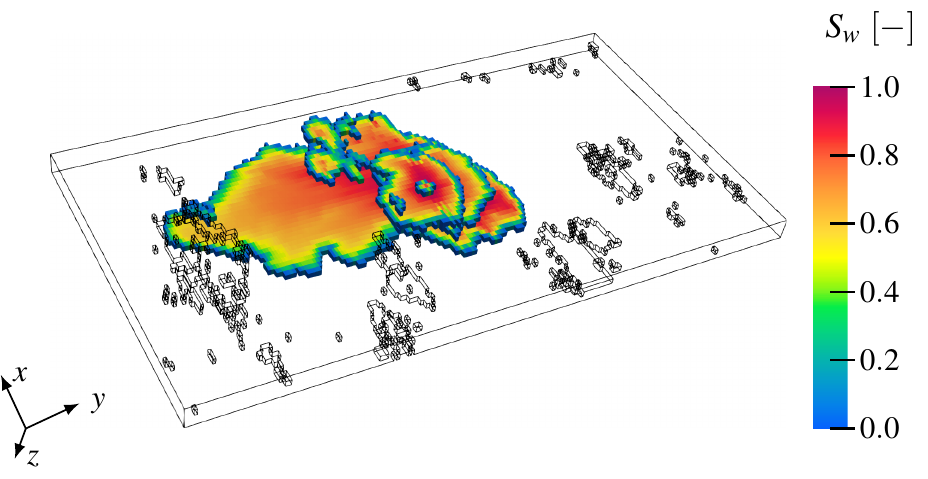}
    }
    \\
    \subfloat[Water velocity]{
    \includegraphics[width=0.47\textwidth]{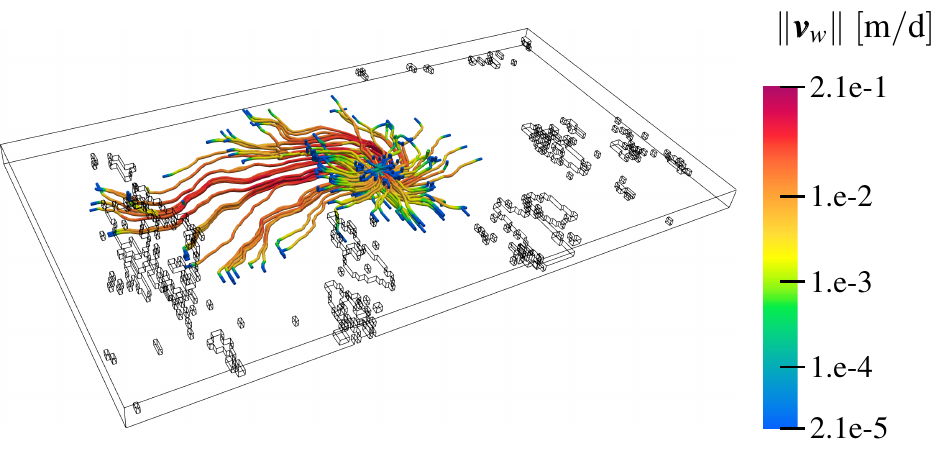}
    }
    \caption{Test 2b: Some model insights at the end of the simulation ($t=1,200$ d).}
    \label{fig:T2b_insights_2_phase}
\end{figure}
This application is the evolution of Test 2a with the introduction of gravity in the model to evaluate its effect on the performance of the BCPR preconditioner. The maximum time step size, equal to 1 day, is smaller than in the previous tests and has been chosen to limit the number of nonlinear iterations to 5 or 6 at most. This results in a larger number of time steps, up to 1,222, although the simulated time interval has been reduced to 1,200 days ($\approx 3.3$ years). The CFL number reaches a value of 24.4 during the simulation. The final outcome from the model is depicted in Figure~\ref{fig:T2b_insights_2_phase}.

\begin{figure*}[tb]
    \centering
    \subfloat[Nonlinear iterations vs.\ time step]{
    \includegraphics[scale=0.9]{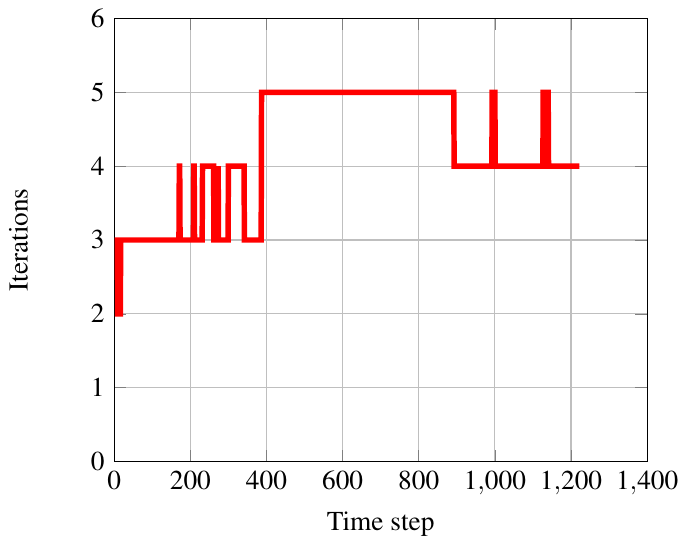}
    } \quad
    \subfloat[Linear iterations vs.\ time step]{
    \includegraphics[scale=0.9]{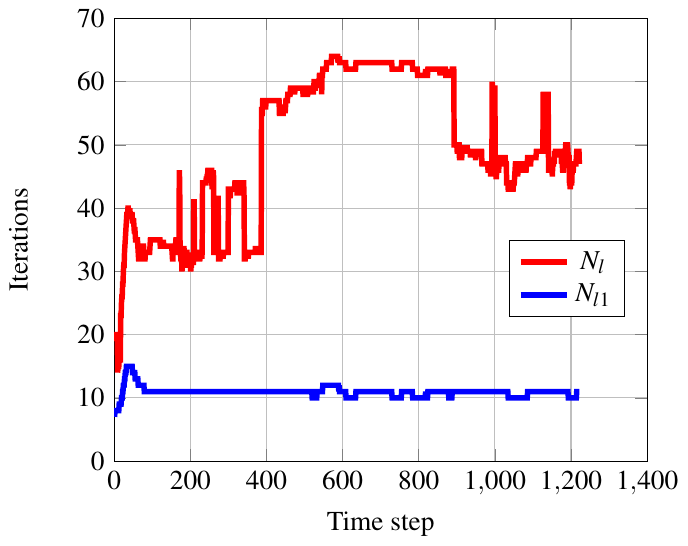}
    } \\
    \subfloat[CPU time vs.\ time step]{
    \includegraphics[scale=0.9]{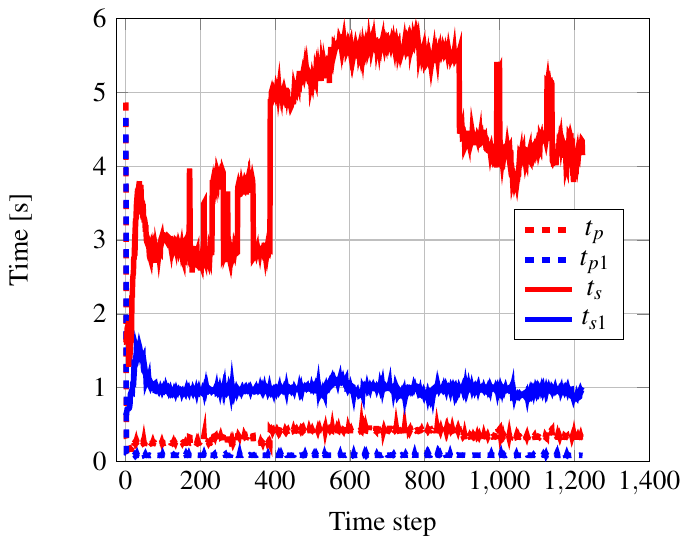}
    }
    \caption{Test 2b: Number of nonlinear (a) and linear (b) iterations and CPU time (c) vs.\ the time step index during the full run 2 of Table~\ref{tab:T2b_static}.}
    \label{fig:T2b_sim_insights}
\end{figure*}

Table~\ref{tab:T2b_static} reports the performance of the BCPR preconditioner with the static approach for building $\widetilde{F}$. As compared to Test 2a, the introduction of gravity does not seem to affect the setup strategy of the preconditioner. Pattern A enables a significant decrease (17.63\%) in the number of linear iterations with respect to the original patch (which allows for convergence in this test) and is again the optimal choice. Further expansions produce a larger cost in terms of CPU time, which, nevertheless, remains close to $1\ \text{s}$ for the solution of a single system after approximately 10 iterations. Although the presence of gravity does not appear to significantly influence the setup and convergence of the linear solver, it does affect that of the nonlinear solver, as shown in Figure~\ref{fig:T2b_sim_insights}. The number of nonlinear iterations, in fact, is larger than in Test 2a, though the maximum time step size is reduced. The other panels confirm that the performance of the BCPR preconditioner is overall stable during the simulation, being the average number of linear iterations per nonlinear step almost constant after the maximum time step size is reached.
\begin{table*}[tb]
    \caption{Test 2b: Numerical performance of the static technique.}
    \centering
    \begin{tabular}{cccccccccccc}
    \toprule
    \# & Pat & $R_S$ & $\hat{N}_N$ & $\hat{N}_l$ & $\hat{t}_p$ [s] & $\hat{t}_s$ [s] & $\hat{t}_t$ [s] & $\overline{N}_{l1}$ & $\overline{t}_{p1}$ [s] & $\overline{t}_{s1}$ [s] & $\overline{t}_{t1}$ [s] \\
    \midrule
    1  & Orig & 1   & 5,105 & 73,148 & 341.0 & 7,064.5 & 7,405.5 & 13.5 & 0.1 & 1.4 & 1.5 \\
    2  &  A   & 1.4 & 5,099 & 60,249 & 433.9 & 5,338.3 & 5,772.2 & 10.9 & 0.1 & 1.0 & 1.1 \\
    3  &  C   & 1.7 & 5,099 & 56,953 & 504.9 & 5,665.6 & 6,170.5 & 10.4 & 0.1 & 1.0 & 1.1 \\
    4  &  E   & 2.0 & 5,099 & 57,413 & 612.0 & 5,823.7 & 6,435.7 & 10.5 & 0.1 & 1.1 & 1.2 \\
    5  &  F   & 2.4 & 5,100 & 57,276 & 759.6 & 5,446.0 & 6,205.6 & 10.4 & 0.1 & 1.0 & 1.1 \\
    \bottomrule
    \end{tabular}
    \label{tab:T2b_static}
\end{table*}

\subsection{Test 3: Dome reservoir with heterogeneous and anisotropic permeability}
\label{sec:test_3}
\begin{figure}
    \centering
    \subfloat[Water saturation]{
    \includegraphics[width=0.47\textwidth]{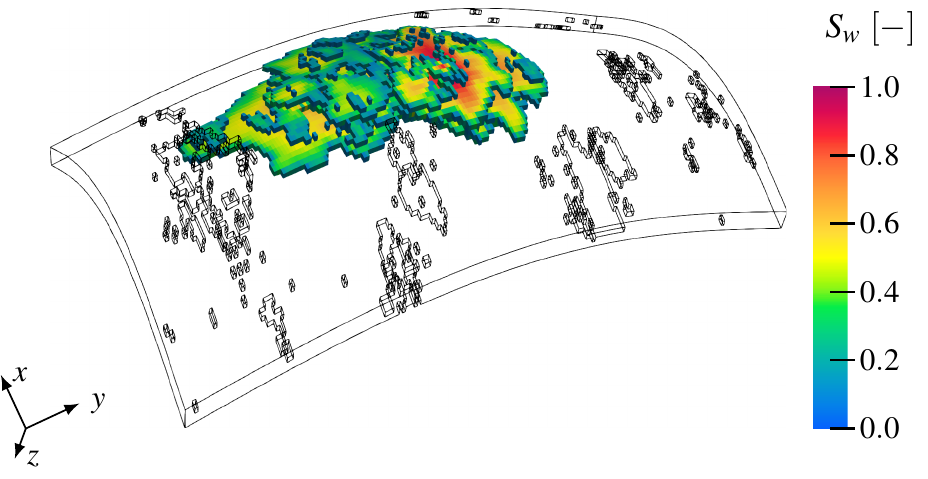}
    }
    \\
    \subfloat[Water velocity]{
    \includegraphics[width=0.47\textwidth]{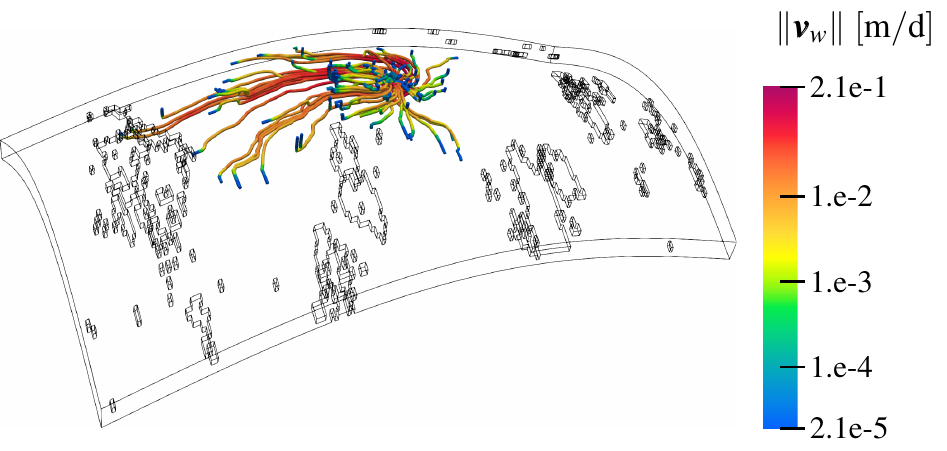}
    }
    \caption{Test 3: Some model insights at the end of the simulation ($t=1,500$ d).}
    \label{fig:T3_insights_2_phase}
\end{figure}

In this application, the simulated temporal domain spans 1,500 days ($\approx 4.1$ years) of continuous production/injection, subdivided in 412 time steps with $\Delta t_{\max}=4$ days and a maximum CFL number equal to 13.2. An overview of the model outcome at the end of the simulation is provided in Figure~\ref{fig:T3_insights_2_phase}, while Tables~\ref{tab:T3_static} and \ref{tab:T3_dynamic} report the results of the numerical investigation, carried out using both the static and dynamic variants for the Schur complement approximation, respectively.
Using a non-Cartesian grid with heterogeneity and anisotropy makes the linearized problem computationally more challenging, as first observed in the comment about Figure~\ref{fig:iter_AGMG_blocks}b and confirmed by the increase in the number of linear iterations in $\overline{N}_{l1}$ (compare Tables~\ref{tab:T3_static} and~\ref{tab:T3_dynamic} to~\ref{tab:T1_static} and~\ref{tab:T2a_static}). The computation and application cost of the preconditioner is also larger than in the previous tests. In fact, blocks $J_{\pi \pi}$ and $\widetilde{S}$ with the original pattern are 3.6 and 1.9 times denser, respectively. This is due to the fact that, with a non-Cartesian grid, the elementary matrix $B^E$ (see equation~\eqref{eq:B_ij_2ph}) is generally full.
This has also an effect on $\widetilde{F}$, for which including the lateral faces in the element-to-face pattern produces now a different approximation (compare runs 2,3 and 4,5 in Table~\ref{tab:T3_static}). 

An appreciable reduction in the total number of linear iterations occurs in run 2 (pattern A), as in the previous tests, however the total CPU time does not decrease correspondingly because of the $R_S$ growth. 
Further pattern expansions turn out to be actually detrimental for the solver speed. Hence, the optimal setup makes use of either the original or the A pattern, with up to 6 additional non-zeros introduced in each column of $\widetilde{F}$.

Based on these results, we might argue that the optimal performance of the dynamic technique is achieved with a similar number of new column entries, $n_{\textup{ent}}$. This is confirmed by a two-step sensitivity analysis in which we tuned $n_{\textup{ent}}$ to find its optimal value, before adjusting $n_{\textup{add}}$, i.e., the number of entries added at each iteration. Table~\ref{tab:T3_dynamic} summarizes the outcome of this investigation. In the first six runs, we fix $n_{\textup{add}}=2$. 
The total number of linear iterations decreases as $n_{\textup{ent}}$ grows up to 6, then the performance deteriorates; therefore, this is the optimal value sought. In the second part of the analysis, where $n_{\text{add}}$ is tuned (runs 7-11), we observe that, as this parameter grows, $\hat{t}_p$ decreases, since fewer iterations are required to compute each column of $\widetilde{F}$. The best results are given for $n_{\textup{add}}=3$ (run 8), then the number of linear iterations starts to increase slowly, reaching the maximum when $n_{\textup{add}}=6$, i.e., all the new column entries are introduced at the first iteration. This result is consistent with the sensitivity analysis performed in~\cite{Nardean2021} on the same non-Cartesian discretization, where we observed that at least two iterations are needed for the decoupling factor approximation to be effective. Overall, the solver performance in run 8 improves the one recorded in Table~\ref{tab:T3_static} both in terms of iteration count and total CPU time.

\begin{table*}[tb]
    \caption{Test 3: Numerical performance of the static technique.}
    \centering
    \begin{tabular}{cccccccccccc}
    \toprule
    \# & Pat & $R_S$ & $\hat{N}_N$ & $\hat{N}_l$ & $\hat{t}_p$ [s] & $\hat{t}_s$ [s] & $\hat{t}_t$ [s] & $\overline{N}_{l1}$ & $\overline{t}_{p1}$ [s] & $\overline{t}_{s1}$ [s] & $\overline{t}_{t1}$ [s] \\
    \midrule
    1  & Orig & 1   & 1,235 & 27,915 & 125.4 & 4,924.7 & 5,050.1 & 24.4 & 0.1 & 4.3 & 4.4 \\
    2  &  A   & 2.0 & 1,235 & 23,879 & 233.6 & 4,797.0 & 5,030.6 & 20.8 & 0.2 & 4.2 & 4.4 \\
    3  &  B   & 2.3 & 1,238 & 28,903 & 266.2 & 7,714.8 & 7,981.0 & 25.2 & 0.2 & 6.8 & 7.0 \\
    4  &  C   & 2.9 & 1,235 & 23,705 & 343.2 & 5,748.9 & 6,092.1 & 20.6 & 0.3 & 5.0 & 5.3 \\
    5  &  D   & 4.3 & 1,224 & 23,268 & 504.7 & 7,671.6 & 8,176.3 & 19.3 & 0.4 & 6.4 & 6.8 \\
    6  &  E   & 3.7 & 1,235 & 23,722 & 456.7 & 7,803.3 & 8,260.0 & 20.5 & 0.4 & 6.8 & 7.2 \\
    7  &  F   & 4.5 & 1,235 & 23,673 & 558.4 & 8,767.7 & 9,326.1 & 20.5 & 0.5 & 7.6 & 8.1 \\
    \bottomrule
    \end{tabular}
    \label{tab:T3_static}
\end{table*}
\begin{table*}[tb]
    \caption{Test 3: Numerical performance of the dynamic technique.}
    \centering
    \begin{tabular}{ccccccccccccc}
    \toprule
    \# & $n_{\textup{ent}}$ & $n_{\textup{add}}$ & $R_S$ & $\hat{N}_N$ & $\hat{N}_l$ & $\hat{t}_p$ [s] & $\hat{t}_s$ [s] & $\hat{t}_t$ [s] & $\overline{N}_{l1}$ & $\overline{t}_{p1}$ [s] & $\overline{t}_{s1}$ [s] & $\overline{t}_{t1}$ [s] \\
    \midrule
    1 & 2  & 2 & 1.4 & 1,235 & 25,420 & 176.7 & 5,191.4 & 5,368.1 & 22.1 & 0.1 & 4.6 & 4.7 \\
    2 & 4  & 2 & 1.7 & 1,235 & 24,429 & 233.0 & 5,308.3 & 5,541.3 & 21.1 & 0.2 & 4.6 & 4.8 \\
    3 & 6  & 2 & 2.0 & 1,235 & 23,151 & 278.4 & 4,409.9 & 4,688.3 & 20.2 & 0.2 & 3.9 & 4.1 \\
    4 & 8  & 2 & 2.2 & 1,235 & 23,881 & 319.7 & 6,043.2 & 6,362.9 & 20.3 & 0.2 & 5.2 & 5.4 \\
    5 & 10 & 2 & 2.5 & 1,235 & 24,682 & 356.8 & 6,869.6 & 7,226.4 & 21.1 & 0.3 & 5.9 & 6.2 \\
    6 & 12 & 2 & 2.7 & 1,235 & 25,878 & 390.1 & 7,577.3 & 7,967.4 & 21.9 & 0.3 & 6.5 & 6.8 \\
    \midrule
    7 & 6  & 1 & 2.0 & 1,235 & 23,299 & 301.1 & 4,493.2 & 4,794.3 & 20.3 & 0.2 & 3.9 & 4.1 \\
    8 & 6  & 3 & 2.0 & 1,235 & 22,904 & 272.7 & 4,356.7 & 4,629.4 & 20.1 & 0.2 & 3.9 & 4.1 \\
    9 & 6  & 4 & 2.0 & 1,235 & 23,179 & 275.2 & 4,471.2 & 4,746.4 & 20.3 & 0.2 & 4.0 & 4.2 \\
    10 & 6 & 5 & 2.0 & 1,235 & 23,910 & 271.6 & 4,683.7 & 4,955.3 & 20.5 & 0.2 & 4.0 & 4.2 \\
    11 & 6 & 6 & 2.0 & 1,235 & 24,404 & 252.1 & 4,903.1 & 5,155.2 & 21.1 & 0.2 & 4.3 & 4.5 \\
    \bottomrule
    \end{tabular}
    \label{tab:T3_dynamic}
\end{table*}

\section{Discussion and conclusions}
The efficient solution of the systems of linearized equations~\eqref{eq:linearized_syst_prot}, originating from an original MHFE-FV discretization of the PDEs that govern the classical two-phase flow model in compressible media, was the main objective of this work. To this end, we designed a preconditioning algorithm built on top of the classical CPR approach and adapted to the $3 \times 3$ block structure of the Jacobian matrix~\eqref{eq:linearized_syst_prot} by introducing a block preconditioner for the $2 \times 2$ pressure subproblem $J_{PP}$. This allows achieving a higher solver robustness by exploiting the inner block structure instead of using a single AMG for the whole $J_{PP}$ part. 
The resulting preconditioner succeeds in incorporating block preconditioning within a global CPR-like algorithm and is thus labeled Block CPR (BCPR). The new core of this tool is the block preconditioner for $J_{PP}$, whose development builds on our previous works~\citep{Nardean2020,Nardean2021} about the Explicit Decoupling Factor Approximation (EDFA) preconditioner. This technique relies on the approximation of the Schur complement with inexact versions of the Jacobian decoupling factors, making use of appropriate restriction/prolongation operators built upon non-zero patterns created either statically or dynamically, as also proposed in~\cite{Ferronato2019,Franceschini2021}. Two additional modifications to the native CPR algorithm have been introduced. Due to the poor performance of a global ILU preconditioner 
and in view of a massive parallelization of the algorithm, the Jacobi preconditioner has been used at the global stage.
The second modification concerns the order of the stages, which has been inverted, following the approaches in~\cite{Bui2017,Roy2020}. The aim is to introduce some kind of decoupling effect of pressure from saturation, since this task is not explicitly performed in a preprocessing step to preserve the original algebraic properties of $J_{PP}$. In essence, we inverted the classical CPR algorithmic structure and modified both preconditioners for the stages. The resulting BCPR algorithm can be cast in an AMG-like framework, where the global stage corresponds to a Jacobi pre-smoothing step, followed by a coarse level solver on pressures carried out by the EDFA tool.

The testing phase in Section~\ref{sec:num_res} showed that the BCPR preconditioner is overall robust and efficient. Its setup consists mainly of tuning the parameters controlling the behavior of the inner EDFA tool, which represents the original core of the BCPR technique, although other elements can be modified such as the AMG function or the global-stage preconditioner. The numerical experiments on different test cases showed that only a few additional entries need to be included in the original column sparsity pattern (Figure~\ref{fig:static_pat}a) to obtain an accurate and sparse, thus efficient, approximation of $F$, hence of the Schur complement itself. These results are consistent with the previous outcome obtained for the EDFA preconditioner alone in a single-phase problem~\cite{Nardean2020,Nardean2021}. Specifically, with a Cartesian grid and regardless of the rock property distribution (Test 1-2), the static technique with pattern A is preferable. The introduction of gravity into the flow model (Test 2b), while causing the nonlinear problem to be computationally more challenging, does not significantly affect the overall behavior and setup of the linear solver. Moving from a Cartesian to a non-Cartesian tessellation modifies the stencil of some blocks in the Jacobian matrix, in particular $J_{\pi \pi}$ and $J_{p \pi}$, and results in a denser approximation of the Schur complement, while using the same pattern, and in higher computation and application cost of the BCPR preconditioner. In this setting, expanding the original pattern with the dynamic technique proved to be more efficient than the static approach, but, once more, only a few additional column entries in $\widetilde{F}$ are required for an optimal performance of the BCPR preconditioner.


This work 
is the basis for future developments, mainly following two principal directions: expanding the underlying physical model and strengthening the implementation.
The BCPR preconditioner has been developed using the two-phase flow problem as a reference, with incompressible fluids and no capillarity effects. However, fluid compressibility and capillarity can be introduced in the flow model as a first step towards an extension of the BCPR framework to multiphase black oil models and, in a farther future, compositional simulations. The good performance of the proposed algorithm also needs to be confirmed for those more sophisticated flow models.
On the other hand, we want to develop a more efficient version of the BCPR preconditioner, derived from the Matlab prototype, using low-level programming languages, such as C++, with the ultimate goal of including it in our in-house built QASR simulator~\cite{Li2021a}. This would also allow linking mainstream AMG libraries to the code and carry out more significant analyses on the CPU time consumption, especially on larger models.
Moreover, ongoing work is being conducted to better understand the benefits of inverting the order of the stages in the BCPR algorithm for its overall effectiveness, by performing a theoretical error analysis accompanied by a thorough numerical investigation.



\section*{Declarations}
\subsection*{Funding}
This publication was supported by the National Priorities Research Program grant NPRP11S-1210-170079 from Qatar National Research Fund.

\subsection*{Competing interests}
The authors have no competing interests to declare that are relevant to the content of this article.

\bibliography{Bib_ECMOR_2022}

\begin{thebibliography}{10}
\providecommand{\doi}[1]{\url{https://doi.org/#1}}
\bibcommenthead

\bibitem{Magras2001}
Magras JF, Quandalle P, Bia P.
\newblock {High-performance reservoir simulation with parallel ATHOS}.
\newblock In: SPE Reserv. Simul. Symp. Houston, Texas, USA: Society of
  Petroleum Engineers; 2001. p. SPE--66342--MS.
\newblock Available from:
  \url{https://onepetro.org/spersc/proceedings/01RSS/All-01RSS/Houston,
  Texas/133525}.

\bibitem{Hu2013}
Hu X, Wu S, Wu XH, Xu J, Zhang CS, Zhang S, et~al.
\newblock {Combined preconditioning with applications in reservoir simulation}.
\newblock Multiscale Model Simul. 2013 jan;11(2):507--521.
\newblock \doi{10.1137/120885188}.

\bibitem{Esler2021}
Esler K, Gandham R, Patacchini L, Garipov T, Samardzic A, Panfili P, et~al.
\newblock {A graphics processing unit–based, industrial grade compositional
  reservoir simulator}.
\newblock SPE J. 2021 oct;p. SPE--203929--PA.
\newblock \doi{10.2118/203929-PA}.

\bibitem{Wallis1983}
Wallis JR.
\newblock {Incomplete Gaussian elimination as a preconditioning for generalized
  conjugate gradient acceleration}.
\newblock In: SPE Reserv. Simul. Symp. San Francisco, California: Society of
  Petroleum Engineers; 1983. p. 325--334.
\newblock Available from: \url{http://www.onepetro.org/doi/10.2118/12265-MS}.

\bibitem{Wallis1985}
Wallis JR, Kendall RP, Little TE.
\newblock {Constrained residual acceleration of conjugate residual methods}.
\newblock In: SPE Reserv. Simul. Symp. Dallas, Texas: Society of Petroleum
  Engineers; 1985. p. SPE--13536--MS.
\newblock Available from: \url{http://www.onepetro.org/doi/10.2118/13536-MS}.

\bibitem{Zhou2013}
Zhou Y, Jiang Y, Tchelepi HA.
\newblock {A scalable multistage linear solver for reservoir models with
  multisegment wells}.
\newblock Comput Geosci. 2013 apr;17(2):197--216.
\newblock \doi{10.1007/s10596-012-9324-0}.

\bibitem{Garipov2018}
Garipov TT, Tomin P, Rin R, Voskov DV, Tchelepi HA.
\newblock {Unified thermo-compositional-mechanical framework for reservoir
  simulation}.
\newblock Comput Geosci. 2018 aug;22:1039--1057.
\newblock \doi{10.1007/s10596-018-9737-5}.

\bibitem{Cremon2020}
Cremon MA, Castelletto N, White JA.
\newblock {Multi-stage preconditioners for thermal–compositional–reactive
  flow in porous media}.
\newblock J Comput Phys. 2020 oct;418:109607.
\newblock \doi{10.1016/j.jcp.2020.109607}.

\bibitem{Klemetsdal2020}
Klemetsdal {\O}S, M{\o}yner O, Lie KA.
\newblock {Accelerating multiscale simulation of complex geomodels by use of
  dynamically adapted basis functions}.
\newblock Comput Geosci. 2020 apr;24(2):459--476.
\newblock \doi{10.1007/s10596-019-9827-z}.

\bibitem{Lie2019}
Lie KA.
\newblock {An introduction to reservoir simulation using MATLAB/GNU Octave}.
\newblock Cambridge, United Kingdom: Cambridge University Press; 2019.
\newblock Available from:
  \url{https://www.cambridge.org/core/product/identifier/9781108591416/type/book}.

\bibitem{Alvestad2022}
Alvestad J, Baxendale D, Bao K, Blatt M, Hove J, Lauser A, et~al.
\newblock {OPM flow: Reference manual}.
\newblock Oslo, Norway: Equinor ASA; 2022.
\newblock Available from:
  \url{https://opm-project.org/wp-content/uploads/2022/05/OPM{\_}Flow{\_}Reference{\_}Manual{\_}2022-04{\_}Rev-0{\_}Reduced.pdf}.

\bibitem{Rasmussen2021}
Rasmussen AF, Sandve TH, Bao K, Lauser A, Hove J, Skaflestad B, et~al.
\newblock {The open porous media flow reservoir simulator}.
\newblock Comput Math with Appl. 2021 jan;81:159--185.
\newblock \doi{10.1016/j.camwa.2020.05.014}.

\bibitem{Schlumberger2020a}
Schlumberger.
\newblock {Eclipse: Technical description}; 2020.

\bibitem{Schlumberger2020}
Schlumberger.
\newblock {Intersect: Technical description}; 2020.

\bibitem{Halliburton2014}
Halliburton.
\newblock {Nexus: Technical reference guide}; 2014.

\bibitem{Lacroix2001}
Lacroix S, Vassilevski YV, Wheeler MF.
\newblock {Decoupling preconditioners in the implicit parallel accurate
  reservoir simulator (IPARS)}.
\newblock Numer Linear Algebr with Appl. 2001 dec;8(8):537--549.
\newblock \doi{10.1002/nla.264}.

\bibitem{Singh2018}
Singh G, Pencheva G, Wheeler MF.
\newblock {An approximate Jacobian nonlinear solver for multiphase flow and
  transport}.
\newblock J Comput Phys. 2018 dec;375:337--351.
\newblock \doi{10.1016/j.jcp.2018.08.043}.

\bibitem{Lacroix2003}
Lacroix S, Vassilevski Y, Wheeler J, Wheeler M.
\newblock {Iterative solution methods for modeling multiphase flow in porous
  media fully implicitly}.
\newblock SIAM J Sci Comput. 2003 jan;25(3):905--926.
\newblock \doi{10.1137/S106482750240443X}.

\bibitem{Cao2005}
Cao H, Tchelepi HA, Wallis JR, Yardumian HE.
\newblock {Parallel scalable unstructured CPR-type linear solver for reservoir
  simulation}.
\newblock In: SPE Annu. Tech. Conf. Exhib. Dallas, Texas: Society of Petroleum
  Engineers; 2005. p. SPE--96809--MS.
\newblock Available from: \url{http://www.onepetro.org/doi/10.2118/96809-MS}.

\bibitem{Gries2014}
Gries S, St{\"{u}}ben K, Brown GL, Chen D, Collins DA.
\newblock {Preconditioning for efficiently applying algebraic multigrid in
  fully implicit reservoir simulations}.
\newblock SPE J. 2014 aug;19(4):726--736.
\newblock \doi{10.2118/163608-PA}.

\bibitem{Nardean2022}
Nardean S, Ferronato M, Abushaikha A.
\newblock {Linear solvers for reservoir simulation problems: An overview and
  recent developments}.
\newblock Arch Comput Methods Eng. 2022 oct;29(6):4341--4378.
\newblock \doi{10.1007/s11831-022-09739-2}.

\bibitem{Roy2020}
Roy T, J{\"{o}}nsth{\"{o}}vel TB, Lemon C, Wathen AJ.
\newblock {A constrained pressure-temperature residual (CPTR) method for
  non-isothermal multiphase flow in porous media}.
\newblock SIAM J Sci Comput. 2020 jan;42(4):B1014--B1040.
\newblock \doi{10.1137/19M1292023}.

\bibitem{White2019}
White JA, Castelletto N, Klevtsov S, Bui QM, Osei-Kuffuor D, Tchelepi HA.
\newblock {A two-stage preconditioner for multiphase poromechanics in reservoir
  simulation}.
\newblock Comput Methods Appl Mech Eng. 2019 dec;357:112575.
\newblock \doi{10.1016/j.cma.2019.112575}.

\bibitem{T.Camargo2021}
{T  Camargo} J, White JA, Castelletto N, Borja RI.
\newblock {Preconditioners for multiphase poromechanics with strong
  capillarity}.
\newblock Int J Numer Anal Methods Geomech. 2021 jun;45(9):1141--1168.
\newblock \doi{10.1002/nag.3192}.

\bibitem{Brezzi1991}
Brezzi F, Fortin M.
\newblock {Mixed and hybrid finite element methods}. vol.~15 of Springer Series
  in Computational Mathematics.
\newblock New York, NY: Springer-Verlag New York; 1991.
\newblock Available from:
  \url{http://link.springer.com/10.1007/978-1-4612-3172-1}.

\bibitem{Abushaikha2017}
Abushaikha AS, Voskov DV, Tchelepi HA.
\newblock {Fully implicit mixed-hybrid finite-element discretization for
  general purpose subsurface reservoir simulation}.
\newblock J Comput Phys. 2017;346:514--538.
\newblock \doi{10.1016/j.jcp.2017.06.034}.

\bibitem{Abushaikha2020}
Abushaikha AS, Terekhov KM.
\newblock {A fully implicit mimetic finite difference scheme for general
  purpose subsurface reservoir simulation with full tensor permeability}.
\newblock J Comput Phys. 2020 apr;406:109194.
\newblock \doi{10.1016/j.jcp.2019.109194}.

\bibitem{Li2021a}
Li L, Abushaikha A.
\newblock {A fully-implicit parallel framework for complex reservoir simulation
  with mimetic finite difference discretization and operator-based
  linearization}.
\newblock Comput Geosci. 2021 oct;\doi{10.1007/s10596-021-10096-5}.

\bibitem{Kuznetsov2003}
Kuznetsov YA.
\newblock {Spectrally equivalent preconditioners for mixed hybrid
  discretizations of diffusion equations on distorted meshes}.
\newblock J Numer Math. 2003 mar;11(1):61--74.
\newblock \doi{10.1163/156939503322004891}.

\bibitem{Maryska2000}
Maryska J, Rozlozn{\'{i}}k M, Tuma M.
\newblock {Schur complement systems in the mixed-hybrid finite element
  approximation of the potential fluid flow problem}.
\newblock SIAM J Sci Comput. 2000 jan;22(2):704--723.
\newblock \doi{10.1137/S1064827598339608}.

\bibitem{Nardean2020}
Nardean S, Ferronato M, Abushaikha AS.
\newblock {A novel and efficient preconditioner for solving Lagrange
  multipliers-based discretization schemes for reservoir simulations}.
\newblock In: ECMOR XVII - 17th Eur. Conf. Math. Oil Recover. Edinburgh:
  European Association of Geoscientists {\&} Engineers; 2020. p. 1--12.
\newblock Available from:
  \url{https://www.earthdoc.org/content/papers/10.3997/2214-4609.202035072}.

\bibitem{Nardean2021}
Nardean S, Ferronato M, Abushaikha AS.
\newblock {A novel block non-symmetric preconditioner for mixed-hybrid
  finite-element-based Darcy flow simulations}.
\newblock J Comput Phys. 2021 jun;p. 110513.
\newblock \doi{10.1016/j.jcp.2021.110513}.

\bibitem{Nardean2022a}
Nardean S, Ferronato M, Abushaikha A.
\newblock {A blended CPR/block preconditioning approach for mixed
  discretization schemes in reservoir modeling}.
\newblock In: ECMOR 2022 Eur. Conf. Math. Geol. Reserv. The Hague, The
  Netherlands: European Association of Geoscientists {\&} Engineers; 2022. p.
  1--15.
\newblock Available from:
  \url{https://www.earthdoc.org/content/papers/10.3997/2214-4609.202244067}.

\bibitem{Coats1980}
Coats KH.
\newblock {An equation of state compositional model}.
\newblock SPE J. 1980 oct;20(5):363--376.
\newblock \doi{10.2118/8284-PA}.

\bibitem{Peaceman1978}
Peaceman DW.
\newblock {Interpretation of well-block pressures in numerical reservoir
  simulation}.
\newblock SPE J. 1978 jun;18(3):SPE--6893--PA.
\newblock \doi{10.2118/6893-PA}.

\bibitem{Chen2006}
Chen Z, Huan G, Ma Y.
\newblock {Computational methods for multiphase flows in porous media}.
\newblock Philadelphia, PA, USA: Society for Industrial and Applied
  Mathematics; 2006.
\newblock Available from:
  \url{http://epubs.siam.org/doi/book/10.1137/1.9780898718942}.

\bibitem{Brooks1964}
Brooks RH, Corey AT.
\newblock {Hydraulic properties of porous media}.
\newblock Fort Collins, Colorado, USA: Colorado State University; 1964.

\bibitem{Aziz1979}
Aziz K, Settari A.
\newblock {Petroleum reservoir simulation}.
\newblock London, United Kingdom: Applied Science Publishers; 1979.

\bibitem{Raviart1977}
Raviart PA, Thomas JM.
\newblock {A mixed finite element method for 2-nd order elliptic problems}.
\newblock In: Galligani I, Magenes E, editors. Math. Asp. Finite Elem. Methods.
  Lect. Notes Math. Berlin, Heidelberg: Springer; 1977. p. 292--315.
\newblock Available from: \url{http://link.springer.com/10.1007/BFb0064470}.

\bibitem{Zhang2021}
Zhang N, Abushaikha AS.
\newblock {An implementation of mimetic finite difference method for fractured
  reservoirs using a fully implicit approach and discrete fracture models}.
\newblock J Comput Phys. 2021 aug;p. 110665.
\newblock \doi{10.1016/j.jcp.2021.110665}.

\bibitem{Maryska1995}
Mary{\v{s}}ka J, Rozlo{\v{z}}n{\'{i}}k M, Tůma M.
\newblock {Mixed-hybrid finite element approximation of the potential fluid
  flow problem}.
\newblock J Comput Appl Math. 1995 nov;63(1-3):383--392.
\newblock \doi{10.1016/0377-0427(95)00066-6}.

\bibitem{Mose1994}
Mos{\'{e}} R, Siegel P, Ackerer P, Chavent G.
\newblock {Application of the mixed hybrid finite element approximation in a
  groundwater flow model: Luxury or necessity?}
\newblock Water Resour Res. 1994 nov;30(11):3001--3012.
\newblock \doi{10.1029/94WR01786}.

\bibitem{Younes2010}
Younes A, Ackerer P, Delay F.
\newblock {Mixed finite elements for solving 2-D diffusion-type equations}.
\newblock Rev Geophys. 2010 mar;48(1):RG1004.
\newblock \doi{10.1029/2008RG000277}.

\bibitem{Younis2011}
Younis RM.
\newblock {Modern advances in software and solution algorithms for reservoir
  simulation} [PhD dissertation].
\newblock Stanford University; 2011.
\newblock Available from:
  \url{https://stacks.stanford.edu/file/druid:fb287kz3299/RMY{\_}PHD{\_}THESIS-augmented.pdf}.

\bibitem{Bui2017}
Bui QM, Elman HC, Moulton JD.
\newblock {Algebraic multigrid preconditioners for multiphase flow in porous
  media}.
\newblock SIAM J Sci Comput. 2017 jan;39(5):S662--S680.
\newblock \doi{10.1137/16M1082652}.

\bibitem{Napov2012}
Napov A, Notay Y.
\newblock {An algebraic multigrid method with guaranteed convergence rate}.
\newblock SIAM J Sci Comput. 2012 jan;34(2):A1079--A1109.
\newblock \doi{10.1137/100818509}.

\bibitem{Notay2010}
Notay Y.
\newblock {An aggregation-based algebraic multigrid method}.
\newblock Electron Trans Numer Anal. 2010;37:123----146.

\bibitem{Notay2012}
Notay Y.
\newblock {Aggregation-based algebraic multigrid for convection-diffusion
  equations}.
\newblock SIAM J Sci Comput. 2012 jan;34(4):A2288--A2316.
\newblock \doi{10.1137/110835347}.

\bibitem{Eisenstat1983}
Eisenstat SC, Elman HC, Schultz MH.
\newblock {Variational iterative methods for nonsymmetric systems of linear
  equations}.
\newblock SIAM J Numer Anal. 1983 apr;20(2):345--357.
\newblock \doi{10.1137/0720023}.

\bibitem{Jiranek2009}
Jir{\'{a}}nek P, Rozlo{\v{z}}n{\'{i}}k M, Gutknecht MH.
\newblock {How to make simpler GMRES and GCR more stable}.
\newblock SIAM J Matrix Anal Appl. 2009 jan;30(4):1483--1499.
\newblock \doi{10.1137/070707373}.

\bibitem{Falgout2002}
Falgout RD, Yang UM.
\newblock {HYPRE: A library of high performance preconditioners}.
\newblock In: Sloot PMA, Hoekstra AG, Tan CJK, Dongarra JJ, editors. Comput.
  Sci. — ICCS 2002. Berlin, Heidelberg: Springer; 2002. p. 632--641.
\newblock Available from:
  \url{http://link.springer.com/10.1007/3-540-47789-6{\_}66}.

\bibitem{Balay2022}
Balay S, Abhyankar S, Adams MF, Benson S, Brown J, Brune P, et~al.
\newblock {PETSc/TAO users manual - ANL-21/39 - Revision 3.17}.
\newblock Argonne National Laboratory; 2022.

\bibitem{Balay1997}
Balay S, Gropp WD, McInnes LC, Smith BF.
\newblock {Efficient management of parallelism in object-oriented numerical
  software libraries}.
\newblock In: Mod. Softw. tools Sci. Comput. Boston, MA: Birkh{\"{a}}user
  Boston; 1997. p. 163--202.
\newblock Available from:
  \url{http://link.springer.com/10.1007/978-1-4612-1986-6{\_}8}.

\bibitem{Christie2001}
Christie MA, Blunt MJ.
\newblock {Tenth SPE comparative solution project: A comparison of upscaling
  techniques}.
\newblock In: SPE Reserv. Simul. Symp. Houston, Texas: Society of Petroleum
  Engineers; 2001. p. 308--317.
\newblock Available from: \url{http://www.onepetro.org/doi/10.2118/66599-MS}.

\bibitem{Saad1986}
Saad Y, Schultz MH.
\newblock {GMRES: A generalized minimal residual algorithm for solving
  nonsymmetric linear systems}.
\newblock SIAM J Sci Stat Comput. 1986 jul;7(3):856--869.
\newblock \doi{10.1137/0907058}.

\bibitem{Saad2003}
Saad Y.
\newblock {Iterative methods for sparse linear systems}.
\newblock Philadelphia, USA: Society for Industrial and Applied Mathematics;
  2003.
\newblock Available from:
  \url{http://epubs.siam.org/doi/book/10.1137/1.9780898718003}.

\bibitem{Ferronato2019}
Ferronato M, Franceschini A, Janna C, Castelletto N, Tchelepi HA.
\newblock {A general preconditioning framework for coupled multiphysics
  problems with application to contact- and poro-mechanics}.
\newblock J Comput Phys. 2019 dec;398:108887.
\newblock \doi{10.1016/j.jcp.2019.108887}.

\bibitem{Franceschini2021}
Franceschini A, Castelletto N, Ferronato M.
\newblock {Approximate inverse-based block preconditioners in poroelasticity}.
\newblock Comput Geosci. 2021 apr;25(2):701--714.
\newblock \doi{10.1007/s10596-020-09981-2}.

\end{thebibliography}


\end{document}